\documentclass[12pt]{article}
\usepackage{amsmath}
\usepackage{amssymb}
\usepackage{amsfonts}
\usepackage{epsfig}
\usepackage[all]{xypic}

\newcommand{\mgn}{\overline{\mathcal{M}}_{g,n}}
\newcommand{\mgno}{\mathcal{M}_{g,n}}

\newcommand{\Mgn}[1]{\overline{\mathcal{M}}_{#1}}
\newcommand{\Mgno}[1]{\mathcal{M}_{#1}}
\newcommand{\hmgnq}[2]{H^{#1} \left( \overline{\mathcal{M}}_{#2}, \mathbb{Q} \right)}
\newcommand{\hmgn}[2]{H^{#1} \left( \overline{\mathcal{M}}_{#2} \right)}
\newcommand{\taut}[2]{T^{#1}_{#2}}
\newcommand{\bgp}[2]{\mathcal{B}^{#1}_{#2}}

\newtheorem{theorem}{Theorem}

\newtheorem{lemma}[theorem]{Lemma}
\newtheorem{proposition}[theorem]{Proposition}

\newtheorem{conjecture}[theorem]{Conjecture}
\newtheorem{definition}[theorem]{Definition}

\title{The fourth tautological group of $\mgn$ \\ and relations with the cohomology}
\author{Marzia Polito}

\begin{document}
\maketitle

\abstract{We give a complete description of the fourth tautological group of the moduli space of pointed 
stable curves, $\mgn$, and prove that for $g \geq 8$ it coincides with the cohomology group with rational 
coefficients. We further give a conjectural upper bound depending on the genus 
for the degree of new tautological relations.}  

\section{Introduction}

Let $\mgn$ be the moduli space of $n$-pointed complex stable algebraic  curves of
genus $g$.

The existence of some degree $4$ relations among tautological classes has  been proved with various methods by
E. Getzler, C. Faber, R. Pandharipande and P.Belorousski, while other relations are obtained as a
consequence of the well known ones in degree $2$.

We actually prove that no other relations can arise, and that for genus $g \geq 8$, the 
cohomology group $\hmgnq{4}{g,n}$ coincides with its tautological subgroup. 
The main results of this paper are formally stated in Theorems \ref{free} and \ref{main}.

\noindent It turns out that new relations appear only in genus up to $5$, whereas for higher genus
all possible relations arise only as a consequence of degree $2$ ones. The proof of this fact 
allows us to suggest in Conjecture \ref{conj} an upper bound depending on the genus 
for higher degree new tautological relations.

As for the methods, E. Arbarello and M. Cornalba proposed in 
\cite{AC} new methods for computing the cohomology groups with rational
coefficients of $\mgn$; their strategy is to establish a strict
relation between the cohomology of the moduli space and the one of the
irreducible components of the boundary, which in turn can be expressed
in terms of moduli spaces of curves with lower genus or with lower number 
of marked points. 
With similar arguments, we establish inductive procedures on genus and/or
number of markings to derive constraints among coefficients in possible relations.

We will therefore be able to give the explicit expression of a new relation in
$\hmgn{4}{3,2}$, whose existence was proved by Faber as a consequence of the existence
of a tautological relation on the open part $\Mgno{3,2}$.
Furthermore, we  will exclude the existence of any relation other than 
the known ones.  
 
A description of $\hmgnq{4}{g}$, for $g\geq 12$,
has been given by D. Edidin in \cite{Ed}, and once the tautological group is known, we can 
adapt his argument to prove that for $g \geq 8$, it coincides with the cohomology. For this,
we make use  of the results by Harer (\cite{Ha}), 
Ivanov (\cite{Iv}) and Loojenga (\cite{Lo1}) on the homology of the mapping
class group.

This paper is extracted from my {\it Tesi di Perfezionamento} at the Scuola Normale Superiore, Pisa.
In the present exposition, many of the calculations will be omitted. The interested reader can find them
all in the thesis (\cite{Po}), available upon request from the author.

I wish to thank my advisor, Enrico Arbarello, as well as Gilberto Bini, Maurizio Cornalba, Carel Faber
and Rahul Pandharipande for many extremely useful conversations.

\section{Stable graphs and  tautological classes\label{notdef}}

To every stable curve $C$ of genus $g$, with $P$ as a set of markings, one can
associate a labelled graph $\Gamma $ in the following way:

\begin{enumerate}
\item  draw a vertex $v$ for every irreducible component $C\left( v\right) $
of the normalization $\widetilde{C}$ of $C$, and label it with the genus 
$g\left( v\right) $ of that component,

\item  draw an edge between two vertices $v_1$, $v_2$ (possibly a loop if 
$v_1=v_2$) whenever the normalization map 
$\nu :\widetilde{C}$ $\rightarrow C$
identifies two points lying respectively in $C\left( v_1\right) $ and 
$C\left( v_2\right) $,

\item  draw a half-edge with vertex $v$ whenever there is a marking in 
$\nu \left( C\left( v\right) \right) $, and label it with the marking's name.
We  denote by $P(v)$ the set of these markings.
\end{enumerate}

We call {\bf marked half-edges} the half-edges constructed in $3$. The total set 
of half-edges is the union of the set of marked half-edges with the set consisting
of the halves of the edges constructed in $2$.

\noindent Let $r\left( v\right) $ be the valence of a vertex, namely the number of
half-edges with  vertex $v$. The stability condition translates to: 
$2g\left( v\right) +r\left( v\right) \geq 3$,
for every vertex $v$. 
The genus of a curve corresponding to the graph $\Gamma $ is 
$g\left( \Gamma \right)=\chi \left( \Gamma \right) 
+\sum_vg\left( v\right) $. Observe that the
construction of the graph is only based on the topological type of the curve.

\begin{definition}
A {\bf $P$-marked stable graph of genus $g$} (briefly a $(g,P)$ graph), is a connected
graph with $n=|P|$ marked half-edges, with the following additional
data: 
\\ 1) each vertex $v$ is labelled with an integer $g(v)$, 
\\ 2) the valence $r(v)$ of any vertex satisfies the stability condition 
      $2g(v)+r(v)\geq 3$, 
\\ 3) there is a bijection between marked half-edges and elements in $P$, 
\\ 4) $g=\chi \left( \Gamma \right) +\sum_vg\left( v\right) $.
\end{definition}

\noindent The codimension of a graph is defined as the number of its edges.

Given a $P$-marked stable graph of genus $g$ and codimension $d$, 
with set of vertices $V$, one can
associate to it a closed stratum of codimension $d$ in $\Mgn{g,P}$.
For every vertex $v \in V$, we let $S\left(v\right) $, denote the
set of unmarked half-edges with vertex $v$.

Let 
$\Mgn{\Gamma} :=
\prod_{v \in V} \Mgno{g\left( v\right),P\left(v\right) \cup S\left( v\right) }$ 

\noindent The map 
$$
\xi_\Gamma :\Mgn{\Gamma} \rightarrow \Mgn{g,P}
$$
is called a boundary map, and has the closed stratum
$\Delta_\Gamma =\xi _\Gamma \left( \Mgn{\Gamma}\right) $ as image.

\noindent The notation $\Mgn{\Gamma} $ will be used also when $\Gamma$ 
is disconnected: if $\Gamma =\Gamma _1\sqcup \Gamma_2$, then 
$\Mgn{\Gamma} =\Mgn{\Gamma_1} \times \Mgn{\Gamma_2}$.

Let $\Gamma$ be a $(g,P)$-graph. 
\begin{definition}\label{stablegraph}
The graph $G$
is a \textbf{ $\Gamma$-graph} if it is the disjoint union
of a collection of $(g(v), P(v) \cup S(v) )$-graphs.
\end{definition}

\noindent Look at a $\Gamma$-graph $G$.
Set $G = \sqcup_{v \in V} G_v $. We can define the map
$$ \Mgn{G} = \prod \Mgn{G_v} \stackrel{\zeta_{G}}{\longrightarrow} \Mgn{\Gamma} $$
as $\zeta_{G}= \{ \xi_{G_v} \}_{v \in V}$.

\bigskip
\bigskip

Let $\pi$ be the forgetful map: 
\begin{eqnarray*}
\pi  &:&\Mgn{g,n+1}\rightarrow \Mgn{g,n} \\
\left[ C,p_1,...,p_n,p_{n+1}\right]  &
\rightarrow &\left[C,p_1,...,p_n\right] 
\end{eqnarray*}

\noindent We will also refer to the map $\pi$ as the universal curve, or
the projection map.

\noindent Let $\sigma _1,...,\sigma _n$ be the  $n$ canonical sections of the forgetful map,
and let $D_i$ be the image of $\sigma_i$. Finally, let $\omega _\pi $ be the relative dualizing sheaf of $\pi $. 

We recall the definition of the
basic cohomology classes in $\Mgn{g,P}$ (see \cite{AC2}):

\begin{definition}
\begin{eqnarray*}
\psi _i &=
&\sigma _i^{*}\left( c_1\left( \omega _\pi \right) \right),i=1,...,n \\
\kappa _a &=&
\pi _{*}\left( \left( c_1\left( \omega _\pi  \left( \sum D_j \right) \right) \right)
^{a+1}\right) ,a=0,...3g-3+n
\end{eqnarray*}
\end{definition}

\noindent The class $\psi _i$ can be interpreted as the first Chern class of the
orbifold bundle whose fiber over the point $\left[ C,p_1,...,p_n\right] $ is
the cotangent bundle to the curve $C$ evaluated at the point $p_i$.

\begin{definition}
A {\bf Mumford class} in $\hmgnq{*}{g,P}$ is a polynomial in the classes 
$\psi _i,\kappa _a$. The {\bf Mumford ring} is 
$$
\mathbb{Q}\left[ \psi _1,...,\psi _n,\kappa _1,...,\kappa _{3g-3+n}\right] .
$$
\end{definition}

\noindent The Mumford ring on a product or a disjoint union of moduli spaces is the
tensor product or the direct sum of the Mumford rings.

\noindent It is worth noticing that the following formula (see Formula 1.7 in \cite{AC2}) 
holds:

$$
\kappa_a =\pi_{*} (\psi_{n+1}^{a+1} ).
$$

\begin{definition}
A {\bf Mumford class }in $H^{*}\left( \mathcal{M}_{g,P},\mathbb{Q}\right) $ is
the pull-back under the inclusion 
$$
\Mgno{g,P}\rightarrow \Mgn{g,P}
$$
of a polynomial in the classes $\psi _i,\kappa _a$.
\end{definition}

\bigskip\

\begin{figure}[h]
\begin{center}
\mbox{\epsfig{file=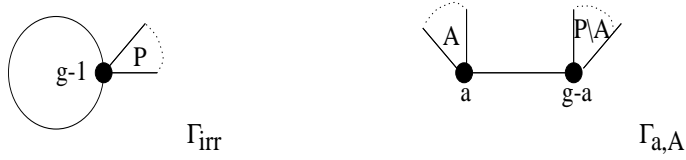,width=9cm,height=2cm}}
\caption{Graphs of codimension $1$\label{codim1}}
\end{center}
\end{figure}

\begin{figure}[h]
\begin{center}
\mbox{\epsfig{file=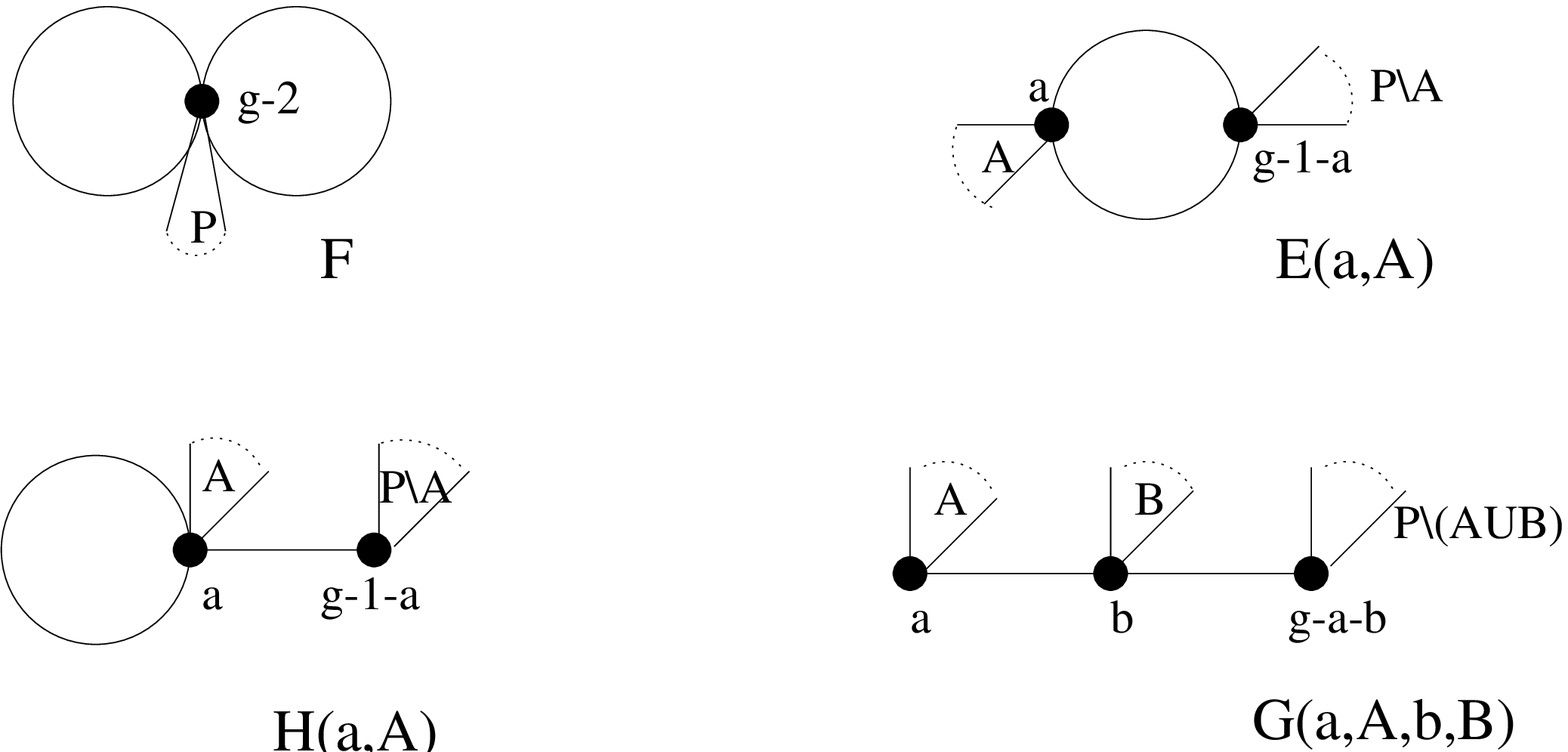,width=10cm,height=6cm}}
\caption{Graphs of codimension $2$\label{codim2}}
\end{center}
\end{figure}

\begin{definition}
A {\bf tautological class} is the push-forward of a Mumford class via a
boundary map. The {\bf $k$-th tautological group} $\taut{k}{g,P}$ is the
subspace of $\hmgnq{k}{g,P}$
generated by these classes.

\end{definition}
In Figures \ref{codim1} and \ref{codim2} we draw all the graphs of 
codimension $1$ and codimension $2$ which we need in our study of
$\taut{4}{g,P}$. In each figure we will also write the name of the
corresponding graph. 
Every time 
half-edges are drawn, one should imagine them labelled with the 
correspondent markings.

\noindent If $p$ is a Mumford class, we use the following notation: 
$$
p|\delta _\Gamma :=\frac{\xi _{\Gamma *}\left( p\right) }{|Aut\Gamma |}. 
$$

\noindent We will often write $\delta_{irr}, \xi_{irr}$ instead of $\delta_{\Gamma_{irr}}, \xi_{\Gamma_{irr}}$,
and $\delta_{a,A}, \xi_{a,A}$ instead of $\delta_{\Gamma_{a,A}}, \xi_{\Gamma_{a,A}}$.

\noindent Degree $4$ autological classes are:

\begin{enumerate}
\item  \textbf{Pure boundary classes}: let $\Gamma $ be a graph of
codimension $2$, then we define: 
$$
\delta _\Gamma :=\frac{\xi _{\Gamma *}\left( 1\right) }{|Aut\Gamma |}
$$

\item  \textbf{Mixed boundary classes:} if codim $\Gamma =1$, and $p$ is a
Mumford class of degree $2$ in $\Mgn{\Gamma} $, then 
$$
p|\delta _\Gamma :=\frac{\xi _{\Gamma *} ( p ) }{|Aut\Gamma |}.
$$

We will often use the following simplified notation:
\begin{itemize}
\item  $\psi _i\delta _{a,A}=\left( \psi _i\otimes 1\right) |\delta _{a,A}
       = \frac{1}{Aut \Gamma_{a,A}} \xi_{a,A*} (\psi_i \otimes 1 )$,

\item  $\psi |\delta _{a,A}=\left( \psi _s\otimes 1\right) |\delta _{a,A}
       = \frac{1}{Aut \Gamma_{a,A}} \xi_{a,A*} (\psi_s \otimes 1 )$,

\item  $\delta _{a,A}|\psi =\left( 1\otimes \psi _t\right) |\delta _{a,A}
       = \frac{1}{Aut \Gamma_{a,A}} \xi_{a,A*} (1 \otimes \psi_t  )
       =\psi |\delta _{g-a,A^c}$,

\item  $\kappa |\delta _{a,A}=\left( \kappa _1\otimes 1\right) |\delta _{a,A}
        = \frac{1}{Aut \Gamma_{a,A}} \xi_{a,A*} (\kappa_1 \otimes 1 )$,

\item $\delta _{a,A}|\kappa = \left( 1\otimes \kappa _1\right) |\delta _{a,A}
        = \frac{1}{Aut \Gamma_{a,A}} \xi_{a,A*} ( 1 \otimes \kappa_1 )
      =\kappa |\delta _{g-,A^c}  $,

\item  $\psi_i \delta _{irr}=\left( \psi _i \right) |\delta _{irr}
        = \frac{1}{Aut \Gamma_{irr}} \xi_{irr*} ( \psi_i)$,

\item  $\psi |\delta _{irr}=\left( \psi _q+\psi _r\right) |\delta _{irr}
        = \frac{1}{Aut \Gamma_{irr}} \xi_{irr*} ( \psi _q+\psi _r )$,

\item  $\kappa_1  \delta _{irr}=\kappa_1 |  \delta _{irr}
       = \frac{1}{Aut \Gamma_{irr}} \xi_{irr*} ( \kappa_1 )$.

\end{itemize}

\item  \textbf{Mumford classes} : these are simply monomials in
Mumford classes (considered as push-forward via the map corresponding
to the trivial graph).
\end{enumerate}

In the mixed boundary classes we intentionally used ambiguous notation. Some
of the classes ($\psi _i\delta _{a,A} , \psi_i \delta _{irr} , \kappa_1  \delta _{irr}$) 
turn out to be written as a product of a codimension $1$
boundary class with a Mumford class. In the proof of the next Proposition 
we will show that the above notation is unambiguous.

\begin{proposition}\label{prodotti}
The image of the map: 
\begin{eqnarray*}
\hmgn{2}{g,P} \times \hmgn{2}{g,P} \rightarrow \hmgn{4}{g,P}  \\
(\alpha ,\beta ) \rightarrow \alpha \cdot \beta 
\end{eqnarray*}
lies in $\taut{4}{g,P}$.

\end{proposition}

\textbf{Proof.} Recall that $\hmgn{2}{g,P}=\taut{2}{g,P}$. 
Two irreducible codimension $1$ boundary classes either coincide
or intersect transversally. In the latter case, it is trivial to check that
their intersection is a linear combination of  tautological pure boundary classes.
 The product of two  Mumford classes is clearly a Mumford class.

\noindent Finally, using the push-pull formula, one is able to express
the product of a Mumford class and a boundary class, and the square
of a boundary class, as linear combination of tautological classes:

\begin{eqnarray*}
\psi _i \cdot \delta _{a,A}=\psi _i | \delta _{a,A} &
\psi _i \cdot \delta _{irr}=\psi _i | \delta _{irr} \\
\kappa _1 \cdot \delta _{a,A}=\kappa _1 | \delta _{a,A}+ \delta_{a,A} | \kappa_1 &
\kappa _1 \cdot \delta _{irr}= \kappa _1 | \delta _{irr}\\
\delta _{a,A}^2  =-\psi | \delta_{a,A} - \delta_{a,A} |\psi & +\left\{ \begin{array}{ll}
\frac 2{|Aut\Gamma _{a,A}|}\delta _{G(g-a,\emptyset ,2a-g,P)} & ifA=P \\ 
\frac 2{|Aut\Gamma _{a,A}|}\delta _{G(a,\emptyset ,g-2a,P)} & ifA=\emptyset
\end{array} \right.  \\
\delta _{irr}^2 & =-\frac 12\xi _{irr*}(\psi _q+\psi _r)+2\delta _F+2\sum
\delta _{E(a,A)}\\
&-\psi |\delta_{irr}+2\delta _F+2\sum \delta _{E(a,A)} 
\end{eqnarray*}

We compute explicitely  one sample case.
Since
\begin{eqnarray*}
\xi _{irr}^{*}(\delta _{irr}) &=&\delta _{irr}
+\sum \delta _{a,A\cup \left\{q\right\} }-\psi _q-\psi _r ,
\end{eqnarray*}
then
\begin{eqnarray*}
2 \delta^2_{irr} = \xi _{irr*}\xi _{irr}^{*}(\delta _{irr}) &=&\frac 12\xi _{irr*}\tilde \xi
_{irr*}(1)+\sum \xi _{irr*}\tilde \xi _{a,A\cup \left\{ q\right\} *}(1)-\xi
_{irr*}(\psi _q+\psi _r),
\end{eqnarray*}
where the symbol $\tilde \xi$ is used for boundary maps of
$\Mgn{g-1, P \cup \{ q,r\} }$. In fact, from now on, when composing two boundary
maps, we will append the second one with the twiddle.

We easily compute: $\frac 12\xi _{irr*}\tilde \xi _{irr*}(1)=\frac 12\xi
_{F*}(1)=4\delta _F$, and then observe that $\xi _{irr}\tilde \xi _{a,A\cup
\left\{ q\right\} }=\xi _{E(a,A)}$ and that the corresponding graph has
automorphism order $2$, unless $P=\emptyset ,a=g/2$, when the order is  $4$.
Moreover, $\xi _{irr*}\tilde \xi _{a,A\cup \left\{ q\right\} *}(1)=\xi
_{irr*}\tilde \xi _{g-a,A^C\cup \left\{ q\right\} *}(1)=|Aut\Gamma
_{E(a,A)}|\delta _{E(a,A)}$. Whenever $|Aut\Gamma _{E(a,A)}|=4$, then by
symmetry only one of the summands above does appear, hence we can write

\begin{eqnarray*}
\delta _{irr}^2 & =-\frac 12\xi _{irr*}(\psi _q+\psi _r)+2\delta _F+2\sum
\delta _{E(a,A)}\\
&-\psi |\delta_{irr}+2\delta _F+2\sum \delta _{E(a,A)}.
\end{eqnarray*}

\begin{flushright}
$\square$
\end{flushright}

\section{Essential tautological classes\label{degclass}}

It is well known that, for genus up to $2$, there are some relations between
degree $2$ tautological classes; thus, certain tautological classes could be 
expressed as linear combination of other ones; they are: 
$\kappa_1$ and $\psi_i , i \in P $  for genera $g=0,1$,
$\kappa_1 $  for genus $g=2$. 

\noindent Moreover, there are Keel's relations among boundary classes in genus $0$.

All these relations reproduce themselves in every genus. The reason is quite clear:
every time there is a relation among tautological classes in the second
cohomology group of a codimension $1$ boundary component, we can push it 
forward to $\hmgn{4}{g,P}$.

In this section we will choose a set of degree $4$ tautological classes
which generate $\taut{4}{g,P}$, by eliminating the above relations.
We will call these classes the {\bf essential} tautological classes.
The set of essential tautological classes will be denoted by $\bgp{4}{g,P}$
and it is obtained from the set of all tautological classes by
removing the {\bf unessential} classes which we are presently going to
list.
 
\noindent The unessential tautological classes are:

\begin{tabular}{llll}
$\psi |\delta _{0,A}=\xi_{0,A*} (\psi_s \otimes 1 )$  & 
$\psi_i \delta _{0,A}=\xi_{0,A*} (\psi_i \otimes 1)$  &
$\kappa |\delta _{0,A} =\xi_{0,A*} (\kappa_1 \otimes 1)$& 
for any $g$, \\
$\psi |\delta _{1,A}= \xi_{1,A*}( \frac{\psi_s \otimes 1 }{|Aut \Gamma_{1,A}|})$ &
$\psi_i \delta _{1,A}= \xi_{1,A*} (\psi_i \otimes 1 )$ &
$\kappa |\delta _{1,A}=\xi_{1,A*}(\frac{\kappa_1 \otimes 1 }{|Aut \Gamma_{1,A}|}) $ &
for any $g$, \\
$\kappa |\delta _{2,A}=\xi_{2,A*}(\frac{\kappa_1 \otimes 1 }{|Aut \Gamma_{2,A}|}) $ &&& 
for any $g$, \\
$\psi | \delta_{irr}= \xi_{irr*}(\frac{\psi_q + \psi_r }{2}) $  &&&for $ g=1,2$\\
$\psi_i  \delta_{irr}= \xi_{irr*}(\frac{\psi_i }{2}) $ &&& for $ g=1,2$, \\
$\kappa | \delta_{irr}= \xi_{irr*}(\frac{1\kappa_1 }{2}) $ &&& for $ g=1,2,3$,\\
$\psi_i^2,\psi_i\psi_j$& $\kappa_1^2, \kappa_1 \psi_i$ && for $g=0,1$,       \\
                       & $\kappa_1^2, \kappa_1 \psi_i$ && for $g=2$.
\end{tabular}

\noindent Moreover, some classes $\delta_{G(0,A,0,B)}$ are unessential (see below); in fact,
in genus $0$ there are Keel's relations (\cite{Ke}) among boundary classes:
we can push them forward by means of the maps 
$$
H^2\left( \overline{\mathcal{M}}_{0,A\cup \left\{ s\right\} }\right) 
\stackrel{\phi _{0,A*}}{\rightarrow }H^4\left( \overline{\mathcal{M}}
_{g,P}\right) 
$$
to obtain the following relations : 
$$
\sum_{\substack{x,y\in B, \\ 
z,w\in C, \\B\cup C=A}}\delta _{G\left( 0,B,0,C\right) }+\delta
_{G\left( 0,C,0,B\right) }=\sum_{\substack{x,z\in B,\\y,w\in C,\\B\cup C=A}}\delta
_{G\left( 0,B,0,C\right) }+\delta _{G\left( 0,C,0,B\right) }, 
$$
$$
\sum_{\substack{x,y\in B,\\z\in C,\\B\cup C=A}}\delta _{G\left( 0,B,0,C\right)}
=\sum_{\substack{x,z\in B,\\ y\in C,\\ B\cup C=A}}\delta _{G\left( 0,B,0,C\right) }. 
$$
We now describe a subset of essential classes of this type; 
if we fix an ordering in $P,$ this induces an
ordering of every subset $A$; a basis for $H^2\left( \overline{\mathcal{M}}
_{0,A\cup \left\{ s\right\} }\right) $ consists of classes $\delta
_{0,\left\{ s\right\} \cup C}$, with $B=A\backslash C$, $|B|\geq 3$, or $
\left| B\right| =2$ and $b<c$ $\forall b\in B,\forall c\in C$. This implies
that we are going to consider only classes $\delta _{G\left( 0,B,0,C\right)} $,
 with $|B|\geq 3$, or $\left| B\right| =2$ and $b<c$ $\forall b\in B$, $\forall c\in C$ .



\section{Pull-back formulas\label{pullbacks}}


In this section we show how to pull back tautological classes to the
codimension $1$ boundary components and to the universal curve.
Let $A$ be a stable $(g,P)$-graph of codimension $1$, as defined 
in the introduction, and let $\Gamma$ be a 
stable connected $(g,P)$-graph of codimension $\leq 2$. 

\noindent We fix our attention on a class of the form 
$ p|\delta_{\Gamma}=\frac{1}{|Aut \Gamma|} \xi_{\Gamma *}(p) .$
We want to describe the boundary components of $\Mgn{A}$ on which the
pull-back $ \xi_{A}^{*} (p|\delta_{\Gamma})$ is supported.

Given any stable $A$-graph $G$,
let $j_{s,t}\left( G\right) $ be the graph
obtained by gluing the half edges $s$ and $t$, and let 
$f_{s,t}\left(G\right) $ be the graph obtained from $j_{s,t}\left( G\right) $ by
collapsing the new edge. 
Via the operation $j_{s,t}$ we are either creating a node on an
irreducible component, or joining two irreducible components at a point.
In either case we are creating a node. Via the operation $f_{s,t}$
we are smoothing the new node.
 
We claim that the boundary components we are looking for correspond to $A$-graphs 
$G$ such that $j_{s,t}\left( G\right) = \Gamma $ or $f_{s,t}\left( G\right)= \Gamma $.
It is very simple to produce graphs $G$ of this sort.

Either  $\Delta_\Gamma \subseteq \Delta_A$, or  $\Delta_\Gamma$ and $\Delta _A$ 
intersect transversally.
If $\Delta_{\Gamma}$ and $\Delta _A$ intersect transversally there must be at least 
a vertex $v$ of $\Gamma$ and a simple Feynman move based at $v$ making $\Gamma$
a degeneration of $A$.
Cutting into a half the edge produced by the Feynman move, and calling the two new
half edges $s$ and $t$, creates a stable $A$-graph $G$ having the property that
$f_{s,t}\left( G\right)= \Gamma$.

Suppose, on the other hand, that $\Delta_\Gamma$ is contained in $\Delta_A$. This
simply means that there is at least one edge of $\Gamma$ cutting which produces 
two half edges $s$ and $t$ and a stable $A$-graph $G$ with the property that
$j_{s,t}\left( G\right)= \Gamma$.
 
\noindent Furthermore we can say that
 $\Delta_{\Gamma} \subseteq \Delta _A$ if and only if there exist a graph $G$
such that  $j_{s,t}\left( G\right)= \Gamma$.

In conclusion, whatever the position of $\Delta_\Gamma$ is with respect 
to $\Delta_A$, we can build a diagram:  

$$\xy
\xymatrix{
\Mgn{G} \ar[r]^{\zeta_{G}} \ar[d]^{\eta_{G}} 
&\Mgn{A}  \ar[d]^{\xi_{A}}\\
\Mgn{\Gamma} \ar[r]^{\xi_{\Gamma}}
& \Mgn{g,P} }
\endxy
$$

for any graph $G$ such that $j_{s,t}\left( G\right)= \Gamma$ or 
$f_{s,t}\left( G\right)= \Gamma$.
The maps $\xi_A$ and $\xi_\Gamma$ are boundary maps, the map
$\zeta_G$ has been defined in section \ref{notdef},
and the map $\eta_G$ consists in joining the two half-edges $s$ and $t$
of the graph $G$.

\noindent Observe that some of these maps could be the identity: e.g if 
$\Gamma = A = \Gamma_{irr}$, then the
trivial $A$ - graph $G$ satisfies
$j_{s,t} ( G )= \Gamma $, and the map $\zeta_{G}$ is the identity.  


\begin{proposition}
\label{formpb} Let $\Gamma $ be any stable graph, of codimension $\leq 2$.
Let $A$ be any graph of codimension $1$. Then the following formula holds: 
$$
\frac{\xi _A^{*}(\xi_{\Gamma* }(p))}{Aut \Gamma} 
 =\sum_{f_{s,t} (G) =\Gamma }  \frac{\zeta_{G*} (\eta_{G}^{*}(p))}{Aut G}
 +\sum_{j_{s,t} (G) =\Gamma } \frac{\zeta_{G*} (\eta_{G}^{*}(p))}{Aut G}
 \cdot  c_1(N_{\xi_A}), 
$$
where we denote by $N_{\xi_A}$ the normal bundle to the map $\xi_A$.

\noindent As usual, we will adopt the simplified notation:
$$
\xi _A^{*}  ( p | \delta_{\Gamma } )
 =\sum_{f_{s,t} (G) =\Gamma } (\eta_{G}^{*} (p)) | \delta_{G}
 +\sum_{j_{s,t} (G) =\Gamma } (\eta_{G}^{*} (p)) | \delta_{G}\cdot  c_1(N_{\xi_A}) .
$$
\end{proposition}


{\bf Proof.} As we already explained, the two cycles $\Delta _\Gamma $ and $\Delta _A$ do not
intersect transversally in $\Mgn{g,P}$ if and only if there
exist a graph $G$ such that $j_{s,t}\left( G\right) =\Gamma $. In this
case, we consider a tubular neighborhood  $T$ of the divisor with
normal crossing $\Delta _A \subset \Mgn{g,P}$.

Consider the diagram:
$$\xy
\xymatrix{
\Mgn{A} \ar[r]^{g_A} \ar[dr]^{\xi_{A}} 
&\Mgn{A} / AutA  \ar[d]^{f_{A}}\\
& \Mgn{g,P} }
\endxy
$$
and the normal bundle  $N_{f_A}$ to the 
map $f_A$. Also observe that $g_A^* N_{f_A} = N_{\xi_A}$.

\noindent Introduce a metric in $N_{f_A}$, construct a tubular neighborhood $\widetilde T$
of its zero section, and extend  $f_A$ in the obvious way to a $C^{\infty}$ map
$$
 \tilde f_A : \widetilde T \longrightarrow T .
$$ 
Take then a sufficiently generic
 $C^{\infty}$ section $s$ of 
 $N_{f_A}$ lying in $\widetilde T$. The composition
$\tilde f_A \circ s \circ  g_A$ yields a 
$C^{\infty}$ map
$$
s_A: \Mgn{A} \longrightarrow \Mgn{g,P}
$$
homotopic to $\xi_A$.

As Poincar\`e duality holds for smooth compact orbifolds, we may pull back 
cycles from $\Mgn{g,P}$ to $\Mgn{A}$.
If $\Delta$ is any irreducible
boundary component, then because of our generic choice of the sections, we have, 
by transverse intersection,

\begin{equation}
s_A^{*} ( [ \Delta ] )  =\sum_i \left[ \Delta_i \right ]   \label{sstar1}
\end{equation}
where the sum ranges over the irreducible components $\Delta_i$ of the preimage
of $\Delta $ in $\Mgn{A}$.

The first step is to describe the irreducible components $\Delta_i$.
We claim that they are of two types, which can combinatorially described as follows.
The first one is simply a cycle $\Delta_G \subset \Mgn{A}$ for each graph $G$
such that $f_{s,t}\left( G\right)= \Gamma$. If $\Delta_A$ and $\Delta_\Gamma$ intersect
transversally, these are the only components $\Delta_i$ appearing in the
above expression. 
If not, the remaining $\Delta_i$'s are all of the form 
$$
\frac{\xi_{G*}\xi_G^* (  c_1(N_{\xi_A}))}{Aut G},
$$
where $G$ is a graph such that $j_{s,t}\left( G\right)= \Gamma$.

Once this is established, we get the Proposition for the case $p=1$, that is:

$$
\xi _A^{*}  (\delta_{\Gamma } )
 =\sum_{f_{s,t} (G) =\Gamma } \delta_{G}
 +\sum_{j_{s,t} (G) =\Gamma } \delta_{G}\cdot  c_1(N_{\xi_A}) .
$$
\begin{figure}
\begin{center}
\mbox{\epsfig{file=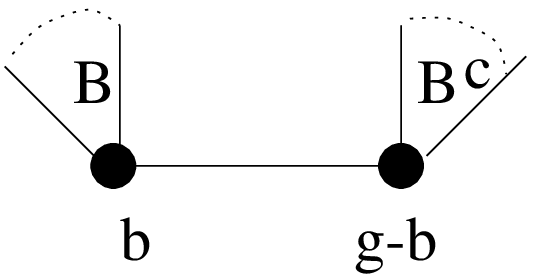,width=3cm,height=1cm}}
\caption{The graph $A=\Gamma=\Gamma_{b,B}$}
\end{center}
\begin{center}
\mbox{\epsfig{file=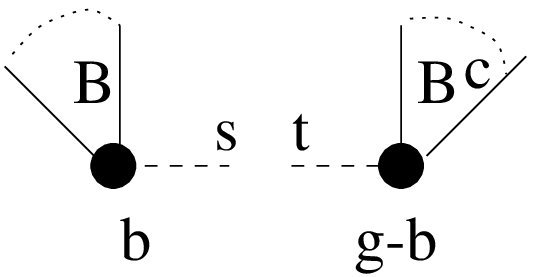,width=3cm,height=1cm}}
\caption{The graph $G$}
\end{center}
\end{figure}
Instead of proving our assertion about the $\Delta_i$'s in general, we 
shall restrict ourselves to some typical examples.
The first example is $\Gamma=A=\Gamma_{b,B}$, with $B \neq \emptyset$, 
$B^c \neq \emptyset$. There is only one $\Delta_i$, which is the zero
locus of a section of the normal bundle to the map $\xi_A$. One may notice
that $\Delta_i$ corresponds to the trivial $A$-graph $G$, drawn on the right, and
that one has that 
$$
 \xi_{b,B}^* (\delta_{b,B})= (\eta_{G}^{*} (1)) | \delta_{G}\cdot  c_1(N_{\xi_A})
                           =c_1(N_{\xi_A}).
$$
This is the standard situation of excess intersection, and there is no 
surprise in finding this term in the general formula of Proposition \ref{formpb}
we are discussing.

\begin{figure}
\begin{center}
\mbox{\epsfig{file=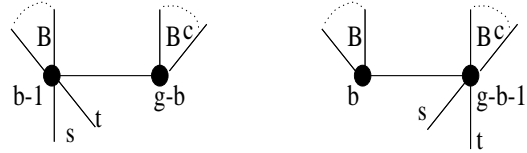,width=7cm,height=2cm}}
\caption{The graphs $G_1$ and $G_2$}
\end{center}
\end{figure}
The opposite situation occurs for example in the formula for 
$$
\xi_{irr}^* (\delta_{b,B})=\xi_{irr}^* (1|\delta_{b,B})
$$
where we further assume that
$b \geq 1, g-b \geq 1$.
There are two components $\Delta_i$, corresponding to the $A$-graphs $G_1$ and 
$G_2$ having the property that 
$f_{s,t}\left( G_i \right)= \Gamma_{b,B}$. In this case
$$
\xi_{irr}^* (\delta_{b,B})= \delta_{G_1} +\delta_{G_2}.
$$
This is the standard situation of transverse intersection.

What is somewhat unexpected in the formula we are discussing, is the mixture
between terms related to excess intersection and terms related to
transverse intersection. To illustrate this phenomenon, let us consider
the case
$$
\xi_{irr}^* (\delta_{F}).
$$ 
The formula in the statement tells us that
$$
\xi _{irr}^{*}\left( \delta _F\right) =-\left( \psi _q+\psi _r\right) \delta_{irr}
+\delta _F
+\sum \delta _{E\left( a,A\cup \left\{ q\right\} \right)}
+\sum \left( \delta _{H\left( a,A\cup \left\{ q\right\} \right) }
+\delta_{H\left( a,A\cup \left\{ r\right\} \right) }\right) ,
$$
where the two sums range over all the possible graphs of the
corresponding type.

\begin{figure}
\begin{center}
\mbox{\epsfig{file=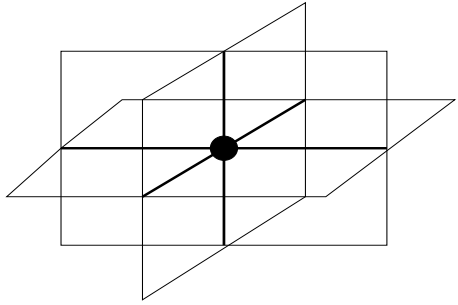,width=6cm,height=4cm}}
\caption{A neighborhood of the three-nodes locus\label{dirr1}}
\end{center}
\end{figure}

The first term is clear: it comes from excess intersection, and corresponds to the
only graph $G$ such that $j_{s,t} (G)=F$, i.e. the graph with one vertex 
of genus $g-2$, one loop, and half-edges with labels in $P \cup \{ s,t \}$.

As a sample case, let us explain the presence of the term $\delta_F$. The presence
of the other terms can be justified by similar arguments.
Draw a picture of 
$\Delta_{irr} $ in a neighborhood of a generic point of the cycle 
$\Delta ^{\prime }$ corresponding to the locus of 
irreducible curves with at least three nodes (Figure \ref{dirr1}). 
We cut it with a codimension three generic subspace, 
in order to draw the picture. The cycle $\Delta ^{\prime }$
is drawn as a triple point of $\Delta _{irr}$, which is locally the union of
three planes, intersecting each other in the three lines belonging to 
$\Delta _F$.

\begin{figure}
\begin{center}
\mbox{\epsfig{file=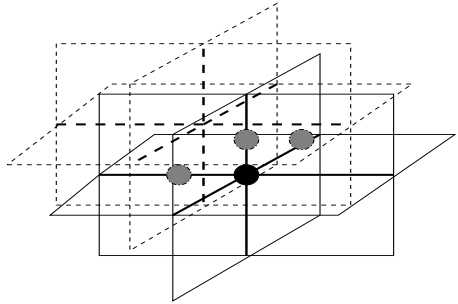,width=6cm,height=4cm}}
\caption{A modified neighborhood of the three-nodes locus \label{dirr2}}
\end{center}
\end{figure}
Now we ``move '' a little bit $\Delta _{irr}$ (Figure \ref{dirr2}),
we call it $\widetilde \Delta_{irr}$, and draw it
with a dotted line. There are three  points of transverse intersection between 
$\widetilde \Delta _{irr}$ and $\Delta _F$. This shows that $s_A^*( \delta_F)$
contains, with multiplicity $1$, the codimension $2$ cocycle in 
$\Mgn{g-1,P \cup \{ s, t \}}$
corresponding to the locus of irreducible two-noded curves, which by abuse of notation is
again denoted by $\Delta_F$.

The formula in the statement, in the case  $p=1$,
$$
\xi _A^{*}  (\delta_{\Gamma } )
 =\sum_{f_{s,t} (G) =\Gamma } \delta_{G}
 +\sum_{j_{s,t} (G) =\Gamma } \delta_{G}\cdot  c_1(N_{\xi_A}) 
$$
is now completely justified.

To prove the general formula we make the following preliminary remark; 
we seek a formula for the pull-back under a $\xi_A$ map of one of the
following classes:

\begin{itemize}
\item  pure boundary classes, hence orbifold Poincar\'e duals of cycles;

\item  $\psi $-mixed classes, hence orbifold Chern classes of bundles
supported on cycles;

\item  $\kappa $-mixed classes. These are linear combinations of the above two
types. In fact, we recall Mumford theorem 
$$
\kappa _1=12\lambda _1+\sum \psi _i-\sum \delta _G,
$$
where the second sum ranges over the set of stable graphs of codimension $1$
, and $\lambda _1$ is the first Chern class of the Hodge bundle; this
implies that $\kappa _1$ is a linear combination of Poincar\'e duals of cycles
and of Chern classes of bundles;

\item  pure Mumford classes, hence polynomials in classes of the above types.
\end{itemize}

\noindent In order to pull-back a tautological class, we first decompose it into  a
linear combination of Mumford classes supported on cycles, and then pull back
each summand separately.

\noindent We therefore seek a formula for 
$$
\xi_A^* (\frac{\xi_{\Gamma *} (c_1 (F))}{Aut \Gamma})
$$
where $F$ is a line bundle on $\Mgn{\Gamma}$. 

\noindent Suppose first that $\Delta_\Gamma$ and $\Delta_A$ intersect transversally. 
Take a sufficiently generic $C^{\infty}$ section $\sigma_F$ of the line bundle
$F$.
For every graph $G$ such that $f_{s,t}(G)=\Gamma$, we denote by $F_G$ the
bundle $\eta_G^* (F)$, and by $\sigma_{F_G}$ its section $\eta_G^*(\sigma_F)$.

By Poincar\'e duality, we can pull back cycles. We claim  that
$$
\xi_A^*(\frac{\xi_{\Gamma *}( [\{ \sigma_F=0 \} ] )}{Aut \Gamma})
= \sum_{f_{s,t}(G)= \Gamma} 
\frac{\xi_G^*( [\{ \sigma_{F_G}=0 \} ] )}{Aut G}.
$$
Let $\Delta$ be a cycle in in $\Mgn{g,P}$ such that 
$$
[\Delta]=\frac{\xi_{\Gamma *}( [\{ \sigma_F=0 \} ] )}{Aut \Gamma};
$$
we can pick
$$
\Delta= \{ x \in \Mgn{g,P}\mid x=\xi_\Gamma (y), \sigma_F (y)=0 \}
$$
with orbifold multiplicity $1$.
Because of transverse intersection of $\Delta_\Gamma$ and $\Delta_A$,
 Formula \ref{sstar1} applies in this case too.
$\Delta$ is a cycle contained in $\Delta_\Gamma$. We therefore  
seek the irreducible components $\Delta_i$ inside the irreducible components
of the preimage of $\Delta_\Gamma$ in $\Mgn{A}$, that is, 
inside the $\Delta_G$'s, where $f_{s,t} (G)= \Gamma$.
One can easily  check that
\begin{eqnarray*}
\Delta_G \cap \xi_A^{-1} (\Delta)& =& \{ z \in \Mgn{A} \cap \Delta_G \mid
       \xi_A (z)=\xi_\Gamma (y) \text{ for some $y$ such that } \sigma_F (y)=0 \} \\
&=&\{ z \in \Mgn{A} \mid z= \zeta_G (w) \text{ for some $w$, } 
\xi_A (z)=\xi_\Gamma (y)\text { for some $y$ such that } \sigma_F (y)=0 \} \\
&=& \{z \in \Mgn{A} \mid z= \zeta_{G}(w) \text{ for some $w$ such that }
\sigma_{F_G} (w)=0 \},
\end{eqnarray*}
again with orbifold multiplicity $1$.
      
Suppose, on the other hand, that $\Delta_\Gamma \subseteq \Delta_A$. 
We need formulas for degree $4$ classes, hence the only new and significant
situation occurs when $\Delta_\Gamma = \Delta_A$, and $\Gamma = A$ is a 
graph of codimension $2$.

From the construction of the map $s_A$, we see that the diagram
$$\xy
\xymatrix{
\Mgn{G} \ar[r]^{\zeta_{G}} \ar[d]^{\eta_{G}} 
&\Mgn{A}  \ar[d]^{s_{A}}\\
\Mgn{\Gamma} \ar[r]^{\xi_{\Gamma}}
& \Mgn{g,P} }
\endxy
$$
commutes only up to homotopy. 
To explain the presence of the transverse intersection terms 
in the pull-back formula, 
$$
\sum_{f_{s,t} (G) =\Gamma }  \frac{\zeta_{G*} (\eta_{G}^{*}(c_1 (F)))}{Aut G},
$$ 
we observe that the induced diagram in cohomology commutes, hence, 
if one chooses suitable sections $\sigma_{F_G}$'s of
the bundles $\eta_G^* (F)$, one can proceed as in the
transverse intersection case.
We now pass to justify the self-intersection term. In our
specific situation this term is
$$
\eta_G^* (c_{1} (F)) \circ c_1(N_{\xi_A}),
$$
in fact, since $\Gamma=A$, the only $A$-graph $G$ such that $j_{s,t}(G)=\Gamma$
is the trivial $A$-graph and  the map $\zeta_G$ is the identity.
The corresponding component in the preimage of $\Delta_\Gamma$ under the map
$s_A$ is the Poincar\'e dual to $c_1(N_{\xi_A})$. Take a section of such
bundle, call it  $\tau$. The component we are looking for is the
Poincar\'e dual of 
$$
\{x \in \Mgn{A} \mid  \sigma_{F_G} (x)=0, \tau (x)=0 \},
$$
that is, the first Chern class of the bundle
$$
\eta_{G}^*(F) \oplus N_{\xi_A},
$$ 
as we claimed.

\begin{flushright}
$\square$
\end{flushright}

\subsection{Formulas for $\pi ^{*}$\label{pistar}}

 
Let 
$$
\pi _A:\overline{{\cal M}}_{g,P\cup A}\rightarrow \overline{{\cal M}} _{g,P} 
$$
be the map forgetting the $A$ markings. We first recall pull-back formulas
for degree $2$ classes (see \cite{AC} and \cite{AC2}).

\smallskip
\begin{tabular}{lll}
$\pi _A^{*}\left( \delta _{c,C}\right) =\sum_{B\subset A}\delta_{c,C\cup B}$ & &
$\pi _A^{*}\left( \psi _i\right) =\psi _i-\sum_{B\subset A}\delta
_{0,B\cup \left\{ i\right\} }$ \\
$\pi _A^{*}\left( \delta _{irr}\right) =\delta _{irr}$ & &
$\pi _A^{*}\left( \kappa _1\right) =\kappa _1-\sum_{i\in A}\psi
_i+\sum_{B\subset A}\delta _{0,B}$
\end{tabular}
\smallskip

\noindent The pull-back formulas for Mumford classes are recursively deduced from 
Formula (1.10) in 
\cite{AC2} and Lemma (1.2) in \cite{AC}; if $\pi :\overline{{\cal M}}
_{0,n}\rightarrow \overline{{\cal M}}_{0,n-1}$ is the forgetful map, then 
\begin{equation}
\psi _i=\pi ^{*}\left( \psi _i\right) +\delta _{0,\left\{ i,n\right\} },
\label{psi}
\end{equation}
and

\begin{equation}
\kappa_i=\pi ^{*}\left( \kappa _i\right) +\psi _n^i.  \label{ki}
\end{equation}

Let us now come to degree $4$ classes.

\noindent Mumford classes are pulled back via formulas \ref{psi} and \ref{ki}:
  
\begin{eqnarray*}
\pi _A^{*}\left( \psi _i^2\right)
&=&\psi _i^2-\sum_{B\subset A}\delta _{0,B\cup \left\{ i\right\} }|\psi +
\text{type }G\text{ classes,}
\end{eqnarray*}

\begin{eqnarray*}
\pi _A^{*}\left( \psi _i\psi _j\right)
&= & \psi _i\psi _j-\psi _j\sum_{B\subset A}\delta _{0,B\cup \left\{ i\right\}}
-\psi _i\sum_{B\subset A}\delta _{0,B\cup \left\{ j\right\} }+\text{type }G
\text{ classes,}
\end{eqnarray*}

\begin{eqnarray*}
\pi _A^{*}\left( \kappa _1\psi _i\right) 
&=&\kappa _1\psi _i-\psi _i\sum_{j\in A}\psi _j-\sum_{B\subset A}
\delta_{0,B\cup \left\{ i\right\} }|\kappa +\sum_{B\subset A,j\in 
A\backslash B}\psi _j\delta _{0,B\cup \left\{ i\right\} }+\sum_{B\subset A}\psi _i\delta
_{0,B} \\
& & +\text{type }G\text{ classes,}
\end{eqnarray*}

\begin{eqnarray*}
\pi _A^{*}\left( \kappa _1^2\right)
&=&\kappa _1^2-2\sum_{i\in A}\kappa _1\psi _i+\sum_{i\in A}\psi
_i^2+2\sum_{i,j\in A,i\neq j}\psi _i\psi _j \\
 & &+2\sum_{B\subset A}\delta _{0,B}|\kappa -2\sum_{B\subset A,i\in
A\backslash B}\psi _i\delta _{0,B}-\sum_{B\subset A}\delta _{0,B}|\psi +
\text{ type }G\text{ classes,}
\end{eqnarray*}

$$
\pi _A^{*}\left( \kappa _2\right) =\kappa _2-\sum_{i\in A}\psi
_i^2+\sum_{B\subset A}\delta _{0,B}|\psi +\text{type }G\text{ classes;} 
$$
this last formula is computed by induction on $|A|$.

With arguments similar to the ones used in Proposition \ref{formpb},
one can easily prove the following:

\begin{proposition}
The following formulas hold: 

\smallskip
\begin{tabular}{lll}
$\pi _A^{*}\left( p|\delta _{irr}\right) =\left( \widetilde{\pi }%
_A^{*}\left( p\right) \right) |\delta _{irr},$ &  & $\pi _A^{*}\left( \delta
_{E\left( c,C\right) }\right) =\sum_{B\subset A}\delta _{E\left( c,C\cup
B\right) },$ \\ 
$\pi _A^{*}\left( p|\delta _{c,C}\right) =\sum_{B\subset A}\left( \widetilde{
\pi }_B^{*}\left( p\right) \right) |\delta _{c,C\cup B}$ &  & $\pi
_A^{*}\left( \delta _{H\left( c,C\right) }\right) =\sum_{B\subset A}\delta
_{H\left( c,C\cup B\right) }$ \\ 
$\pi _A^{*}\left( \delta _F\right) =\delta _F$ &  & $\pi _A^{*}\left( \delta
_{G\left( c,C,d,D\right) }\right) =\sum_{\left( B\cup B^{\prime }\right)
\subset A}\left( \delta _{G\left( c,C\cup B,d,D\cup B^{\prime }\right)
}\right) $
\end{tabular}
\smallskip

where 
\begin{eqnarray*}
\widetilde{\pi }_A &:&\overline{{\cal M}}_{g-1,P\cup A\cup \left\{
q,r\right\} }\rightarrow \overline{{\cal M}}_{g-1,P\cup \left\{ q,r\right\} }
\\
\widetilde{\pi }_B &:&\overline{{\cal M}}_{c,C\cup B\cup \left\{
s\right\} }\times \overline{{\cal M}}_{g-c,\left( P\backslash C\right) \cup
(A\backslash B)\cup \left\{ t\right\} }\rightarrow \overline{{\cal M}}
_{c,C\cup \left\{ s\right\} }\times \overline{{\cal M}}_{g-c,\left(
P\backslash C\right) \cup \left\{ t\right\} }\text{.}
\end{eqnarray*}
\end{proposition}

\begin{flushright}
$\square$
\end{flushright}

\section{Relations in degree $4$\label{newrel}}

New relations arising in degree $4$ appear in $\mgn$ for $ g \leq 5$ 
and for suitable $n$, and can be pulled back with formulas in 
\ref{pistar}.  
They have been computed with different techniques by E. Getzler, R. Pandharipande,
P. Belorousski, and C. Faber. Most of them can be found in the literature, and
we will give below the precise reference.
The existence of some of them follows from \cite{Fa5}, as a consequence of the existence
of tautological relations on $\mgno$, while their explicit expression on $\mgn$
has been recently computed  by C. Faber and privately communicated to the author (\cite{Fa}).
The only exception is the new relation in $\Mgn{3,2}$, whose coefficients will be determined
in section \ref{relat} by the ``pull-back to the boundary'' techniques.

\subsection{Genus $0$}

The only new result is that 
$$
\kappa _2=0\text{ in }H^4\left( \overline{\mathcal{M}}_{0,4}\right) 
$$
for dimension reasons.

\subsection{Genus $1$}

As above, 
$$
\kappa _2=0\text{ in }H^4\left( \overline{\mathcal{M}}_{1,1}\right) . 
$$
Moreover, as observed  by Faber in \cite{Fa3}, 
$$
\delta _{irr}^2=0\text{.} 
$$
There are other relations: the first one originates in 
$H^4\left( \overline{\mathcal{M}}_{1,2}\right) :$ 
$$
\delta _{E\left( 0,\left\{ i\right\} \right) }-\delta _{H\left( 0,\emptyset\right) }=0 , 
$$
as the push-forward of Keel relation with the map $\xi _{irr}:\overline{
\mathcal{M}}_{0,4}\rightarrow \overline{\mathcal{M}}_{1,2}$. The second
one originates in $H^4\left( \overline{\mathcal{M}}_{1,4}\right) $: 
\begin{eqnarray*}
\ \ \ \ 0 &=&12\sum_i\delta _{G\left( 0,\left\{ 1,i\right\} ,1,\emptyset
\right) }-12\sum_i\delta _{G(1,\left\{ i\right\} ,0,\left\{ *\right\}
)}-2\sum_{i,j}\mathbf{\delta }_{G(1,\emptyset ,0,\left\{ i,j\right\} )} \\
&&\ \ \ \ \ \ \ \ \ \ \ \ \ +6\sum_i\delta _{G(1,\emptyset ,0,\left\{
i\right\} )}-2\sum_i\delta _{E(\left\{ 1,i\right\} )}+\sum_i\delta
_{H(\left\{ i\right\} )}+\delta _{H\left( \emptyset \right) }.
\end{eqnarray*}
This was discovered by Getzler (\cite{G1}), while Pandharipande (\cite{Pa})
then proved  it is algebraic.

\subsection{Genus $2$}

Following Mumford (\cite{Mu}), 
$$
60\kappa _2=\delta _F+6\delta _{H\left( 0,\emptyset \right) } 
$$
in $H^4\left( \overline{\mathcal{M}}_{2,0}\right) $.
Faber proves that in $H^4\left( \overline{\mathcal{M}}%
_{2,1}\right) $
$$
\psi _i^2=\frac 1{120}\delta _F+\frac 15\delta _{E\left( 1,\emptyset \right)
}+\frac{13}{120}\delta _{H\left( 0,i\right) }-\frac 1{120}\delta _{H\left(
0,\emptyset \right) }+\frac 75\delta _{G\left( 1,\emptyset ,0,i\right) }. 
$$
Getzler proves in (\cite{G2}) that, in $H^4\left( \overline{\mathcal{M}}%
_{2,2}\right) $, 
\begin{eqnarray*}
\psi _i\psi _j &=&3\psi |\delta _{2,\emptyset }+\frac 1{72}\delta _F+\frac
7{15}\delta _{E\left( 1,\emptyset \right) }+\frac 1{15}\left( \delta
_{E\left( 1,i\right) }+\delta _{E\left( 1,j\right) }\right) \\
&&\ \ \ \ \ +\frac{23}{120}\delta _{H\left( 0,ij\right) }+\frac 1{24}\left(
\delta _{H\left( 0,i\right) }+\delta _{H\left( 0,j\right) }\right) -\frac
1{40}\delta _{H\left( 0,\emptyset \right) }-\frac 1{15}\delta _{H\left(
1,\emptyset \right) } \\
&&\ \ \ \ \ +\frac{13}5\delta _{G\left( 1,\emptyset ,0,ij\right) }+\frac
45\left( \delta _{G\left( 1,i,0,j\right) }+\delta _{G\left( 1,j,0,i\right)
}\right) -\frac 45\delta _{G\left( 0,ij,1,\emptyset \right) }.
\end{eqnarray*}

A new algebraic relation was discovered by Belorousski and Pandharipande (\cite{BP}) 
in $H^4\left( \overline{\mathcal{M}}_{2,3}\right) $: 
\begin{eqnarray*}
\ 0 &=&12\psi |\delta _{2,\emptyset }-6\sum_{i=1}^3\psi |\delta
_{2,i}+6\sum_{i=1}^3\psi _i\delta _{2,i}+\frac 65\delta _{E\left(
1,\emptyset \right) }-\frac 65\sum_{i=1}^3\delta _{E\left( 1,i\right)
}+\frac 25\sum_{i=1}^3\delta _{E\left( 0,i\right) } \\
&&\ +\frac 1{10}\delta _{H\left( 0,123\right) }-\frac
3{10}\sum_{i=1}^3\delta _{H\left( 0,jk\right) }+\frac
3{10}\sum_{i=1}^3\delta _{H\left( 0,i\right) }-\frac 1{10}\delta _{H\left(
0,\emptyset \right) }-\frac 35\delta _{H\left( 1,\emptyset \right) }\ -\frac
15\sum_{i=1}^3\delta _{H\left( 1,i\right) } \\
&&\ -12\delta _{G\left( 2,\emptyset ,0,*\right) }+\frac{12}5\delta _{G\left(
1,\emptyset ,0,123\right) }-\frac{12}5\sum_{i=1}^3\delta _{G\left(
1,i,0,jk\right) }+\frac{24}5\sum_{i=1}^3\delta _{G\left( 1,\emptyset
,0,i\right) } \\
&&-\frac{36}5\sum_{i=1}^3\delta _{G\left( 1,*,0,i\right) }-\frac{36}
5\sum_{i=1}^3\delta _{G\left( 1,\emptyset ,1,\emptyset \right) }+\frac{18}
5\sum_{i=1}^3\delta _{G\left( 1,i,1,\emptyset \right) }-\frac{12}
5\sum_{i=1}^3\delta _{G\left( 1,\emptyset ,1,i\right) }.
\end{eqnarray*}
Here, and from now on, every time we write the symbol $*$ instead
of a marking's name, we mean that any marking which does
not appear elsewhere in the notation could replace the $*$.

\subsection{Genus $3$}

In $H^4\left( \overline{\mathcal{M}}_{3,0}\right) $(\cite{Fa} and \cite{Fa1}): 
\begin{eqnarray*}
\kappa _1^2 &=&-\frac 57\psi |\delta _{irr}-\frac{89}7\psi |\delta
_{2,\emptyset }-\frac 2{35}\delta _F-\frac{94}{35}\delta _{E\left(
1,\emptyset \right) }+\frac{103}{84}\delta _{H\left( 0,\emptyset \right)
}-\frac 27\delta _{H\left( 1,\emptyset \right) }-\frac{22}{35}\delta
_{G\left( 1,\emptyset ,1,\emptyset \right) }, \\
\kappa _2 &=&-\frac 5{42}\psi |\delta _{irr}-\frac{41}{21}\psi |\delta
_{2,\emptyset }+\frac 1{630}\delta _F-\frac{11}{35}\delta _{E\left(
1,\emptyset \right) }+\frac{41}{252}\delta _{H\left( 0,\emptyset \right)
}+\frac 2{105}\delta _{H\left( 1,\emptyset \right) }+\frac 8{35}\delta
_{G\left( 1,\emptyset ,1,\emptyset \right) },
\end{eqnarray*}
whereas in $H^4\left( \overline{\mathcal{M}}_{3,1}\right) $ a new relation
involving $\kappa _1\psi _i$ appears, and the three of them could be written as
follows (\cite{Fa}): 
\begin{eqnarray*}
\kappa _1\psi _i &=&-5\psi _i^2-\frac 17\psi _i\delta _{irr}-
\frac 1{42}\psi|\delta _{irr}-\frac 57\psi _i\delta _{2,i}
-\frac{16}{21}\psi |\delta _{2,i}-\frac{40}{21}\psi \delta _{2,\emptyset }
-\frac 1{630}\delta _F\\ &&
+\frac{13}{21}\delta _{E\left( 0,i\right) }-\frac 9{35}\delta _{E\left( 1,i\right) } 
+\frac{61}{252}\delta _{H\left( 0,i\right) }-\frac 2{105}\delta _{H\left(
1,i\right) } +\frac 4{105}\delta _{H\left( 1,\emptyset \right) }+\frac
4{63}\delta _{H\left( 0,\emptyset \right) } \\ && +\frac{16}{35}\delta _{G\left(
1,i,1,\emptyset \right) }+\frac{61}{21}\delta _{G\left( 1,\emptyset
,0,i\right) }-\frac 8{35}\delta _{G\left( 1,\emptyset ,1,i\right) }, \\
\kappa _1^2 &=&-9\psi _i^2-\frac 27\psi _i\delta _{irr}-\frac{16}{21}\psi
|\delta _{irr}-\frac{10}7\psi _i\delta _{2,i}-\frac{299}{21}\psi |\delta
_{2,i}-\frac{347}{21}\psi \delta _{2,\emptyset }-\frac{19}{315}\delta _F \\ 
&& +\frac{83}3\delta _{E\left( 0,i\right) }-\frac{16}5\delta _{E\left(
1,i\right) } 
+\frac{431}{252}\delta _{H\left( 0,i\right) }-\frac{34}{105}\delta
_{H\left( 1,i\right) }-\frac{22}{105}\delta _{H\left( 1,\emptyset \right) }+
\frac{341}{252}\delta _{H\left( 0,\emptyset \right) } \\ 
&& +\frac 27\delta
_{G\left( 1,i,1,\emptyset \right) }+\frac{389}{21}\delta _{G\left(
1,\emptyset ,0,i\right) }-\frac{38}{35}\delta _{G\left( 1,\emptyset
,1,i\right) }, \\
\kappa _2 &=&-\psi _i^2-\frac 5{42}\psi |\delta _{irr}-\frac{41}{21}\psi
|\delta _{2,i}-\frac{347}{21}\psi \delta _{2,\emptyset }+\frac 1{630}\delta
_F \\ && +\frac 5{21}\delta _{E\left( 0,i\right) }-\frac{11}{35}\delta _{E\left(
1,i\right) } +\frac{41}{252}\delta _{H\left( 0,i\right) }+\frac 2{105}\delta _{H\left(
1,i\right) }+\frac 2{105}\delta _{H\left( 1,\emptyset \right) }+\frac{41}{252
}\delta _{H\left( 0,\emptyset \right) } \\ && +\frac 8{35}\delta _{G\left(
1,i,1,\emptyset \right) }+\frac{41}{21}\delta _{G\left( 1,\emptyset
,0,i\right) }+\frac 8{35}\delta _{G\left( 1,\emptyset ,1,i\right) }.
\end{eqnarray*}

Finally, in $\hmgn{4}{3,2}$, we have:
\begin{eqnarray*}
0&=&\psi_a^2 +\psi_b^2 - \frac65 \psi_a \psi_b - \kappa | \delta_{3, \emptyset} + 5 \psi |\delta_{3, \emptyset}
    - \frac{40}{21} \psi |\delta_{2, \emptyset} + \frac53 \left( \psi |\delta_{2, a}+\psi |\delta_{2, b} \right) \\
&&  - \frac67 \left( \psi_a \delta_{2, a}+\psi_b \delta_{2, b} \right) - \frac{16}{21} \psi |\delta_{2, ab} 
   + \frac{12}{35} \left( \psi_a \delta_{2, ab}+\psi_b \delta_{2, ab} \right) - \frac1{42} \psi | \delta_{irr} \\
&&    + \frac1{35} \left( \psi_a  \delta_{irr}+ \psi_b  \delta_{irr} \right) 
   - \frac1{630} \delta_F + \frac{13}{21} \delta_{E(2, \emptyset)} 
   - \frac4{15} \left( \delta_{E(2, a)} + \delta_{E(2, b)} \right) \\
&& - \frac9{35} \delta_{E(1, \emptyset)} - \frac{34}{105} \delta_{E(1, a)}
   + \frac17 \delta_{H(2, \emptyset)} - \frac2{105} \delta_{H(1, ab)} + \frac4{105} \delta_{H(1, \emptyset)}\\
&& + \frac1{105} \left( \delta_{H(1, a)} + \delta_{H(1, b)} \right)  + \frac4{63} \delta_{H(0, \emptyset)} 
   + \frac{10}{63} \delta_{H(0, ab)}   - \frac5{36} \left( \delta_{H(0, a)} + \delta_{H(0, b)} \right) \\
&& + \frac{40}{21} \delta_{G(2, \emptyset, 0, ab)} - \delta_{G(2, \emptyset, 1, \emptyset)}
   + \frac{16}{35} \delta_{G(1, ab, 1, \emptyset)}  - \frac8{35} \delta_{G(1, \emptyset, 1, ab)} \\
&& - \frac{5}{3} \left( \delta_{G(2, b, 0, a)} +\delta_{G(2, a, 0, b)} \right)
   - \frac{40}{21} \left( \delta_{G(2, \emptyset , 0, a)} +\delta_{G(2, \emptyset, 0, b)} \right).
\end{eqnarray*}

\subsection{Genus $4$}

In $H^4\left( \overline{\mathcal{M}}_{4,\emptyset }\right) $(\cite{Fa} and 
\cite{Fa2}):

\begin{eqnarray*}
0 &=&\frac{45}2\kappa _1^2-240\kappa _2-7\kappa _1\delta _{irr}+\frac{35}
2\psi |\delta _{irr}-39\kappa |\delta _{3,\emptyset }+\frac{315}2\psi
|\delta _{3,\emptyset }+\frac{45}2\psi |\delta _{2,\emptyset } \\
&&+\delta _F+13\delta _{E\left( 2,\emptyset \right) }-\frac{105}8\delta
_{H\left( 0,\emptyset \right) }+2\delta _{H\left( 1,\emptyset \right)
}+5\delta _{H\left( 2,\emptyset \right) }+24\delta _{G\left( 1,\emptyset
,1,\emptyset \right) }+21\delta _{G\left( 1,\emptyset ,2,\emptyset \right) },
\end{eqnarray*}
and since another relation appears in $H^4\left( \overline{\mathcal{M}}_{4,1}\right) $
(\cite{Fa}), we get there the following two relations: 
\begin{eqnarray*}
0 &=&5\kappa _1^2-30\kappa _2-40\kappa _1\psi _i+245\psi _i^2-\kappa
_1\delta _{irr}+7\psi _i\delta _{irr}-2\kappa |\delta _{3,i}+44\psi _i\delta
_{3,i} \\
&&\ -35\psi |\delta _{3,i}-32\kappa |\delta _{3,\emptyset }+175\psi |\delta
_{3,\emptyset }-30\psi _i\delta _{2,i}+95\psi |\delta _{2,i}-85\psi |\delta
_{2,\emptyset } \\
&&\ \ -36\delta _{E\left( 0,i\right) }+24\delta _{E\left( 1,i\right)
}-12\delta _{E\left( 2,i\right) }+\frac{35}{12}\delta _{H\left( 0,\emptyset
\right) }+\delta _{H\left( 1,\emptyset \right) }+5\delta _{H\left(
2,\emptyset \right) }\\
&&-\frac{175}{12}\delta _{H\left( 0,i\right) }+\delta_{H\left( 1,i\right) }
-\delta _{H\left( 2,i\right) } 
 -18\delta _{G\left( 1,i,1,\emptyset \right) }+28\delta _{G\left(
1,i,2,\emptyset \right) }\\
&&+12\delta _{G\left( 2,i,1,\emptyset \right)
}-175\delta _{G\left( 3,\emptyset ,0,i\right) }-10\delta _{G\left(
2,\emptyset ,0,i\right) }+12\delta _{G\left( 1,\emptyset ,1,i\right)
}-4\delta _{G\left( 1,\emptyset ,2,i\right) } \\
0 &=&\frac{25}2\kappa _1^2-180\kappa _2+35\kappa _1\psi _i-\frac{455}2\psi
_i^2-5\kappa _1\delta _{irr}+\frac{35}2\psi |\delta _{irr}-7\psi _i\delta
_{irr} \\
&&\ \ -35\kappa |\delta _{3,i}-49\psi _i\delta _{3,i}+\frac{455}2\psi
|\delta _{3,i}+25\kappa |\delta _{3,\emptyset }-\frac{385}2\psi |\delta
_{3,\emptyset }+60\psi _i\delta _{2,i}-\frac{335}2\psi |\delta _{2,i}+\frac{
385}2\psi |\delta _{2,\emptyset } \\
&&\ \ \ +\delta _F+37\delta _{E\left( 0,i\right) }-35\delta _{E\left(
1,i\right) }+37\delta _{E\left( 2,i\right) }-\frac{455}{24}\delta _{H\left(
0,\emptyset \right) }-5\delta _{H\left( 2,\emptyset \right) }+\frac{385}{24}
\delta _{H\left( 0,i\right) }+7\delta _{H\left( 2,i\right) } \\
&&\ \ +60\delta _{G\left( 1,i,1,\emptyset \right) }-35\delta _{G\left(
1,i,2,\emptyset \right) }+\frac{385}2\delta _{G\left( 3,\emptyset
,0,i\right) }-25\delta _{G\left( 2,\emptyset ,0,i\right) }+49\delta
_{G\left( 1,\emptyset ,2,i\right) }.
\end{eqnarray*}
\subsection{Genus $5$}

Finally, in $H^4\left( \overline{\mathcal{M}}_{5,0}\right) $(\cite{Fa}):

\begin{eqnarray*}
0 &=&\frac{25}2\kappa _1^2-180\kappa _2-5\kappa _1\delta _{irr}+\frac{35}
2\psi |\delta _{irr}-35\kappa |\delta _{4,\emptyset }+\frac{455}2\psi
|\delta _{4,\emptyset }+25\kappa |\delta _{3,\emptyset }\\
&&-\frac{385}2\psi|\delta _{3,\emptyset }-\frac{385}2\delta _{3,\emptyset }|\psi 
\ +\delta _F+37\delta _{E\left( 1,\emptyset \right) }-35\delta _{E\left(
2,\emptyset \right) }  -\frac{455}{24}\delta _{E\left( 0,\emptyset \right)
}\\ && -5\delta _{H\left( 2,\emptyset \right) }+7\delta _{H\left( 3,\emptyset
\right) }-35\delta _{G\left( 1,\emptyset ,2,\emptyset \right) }+49\delta
_{G\left( 1,\emptyset ,3,\emptyset \right) }+25\delta _{G\left( 2,\emptyset
,1,\emptyset \right) }.
\end{eqnarray*}

\section{Degree $4$ relations in the tautological group\label{relat}}

\begin{theorem}
\label{free}For $g\geq 6$, $\bgp{4}{g,P}$ is a basis for $\taut{4}{g,P}$.
For $2 \leq g\leq 5$, the relations among elements of 
$\bgp{4}{g,P}$ are the ones listed in section \ref{newrel}.
\end{theorem}

We will prove this Theorem by induction on $g$.  
We start with a sketchy exposition of an argument which covers 
the cases $g \geq 6$, once the previous ones are established.
Unfortunately, this argument fails to extend to the low genus cases. 
We will therefore give a second, less direct argument. 
The initial cases require more involved computations,
because of the presence of many relations among tautological classes.
We will work out two sample cases in Lemmas \ref{sei} and \ref{due},
and recover the coefficients of the new relation in $\Mgn{3,2}$ in
Proposition \ref{rel32}.


\begin{proposition}
\label{dim1} Suppose that Theorem \ref{free} holds for $g=5$. Then it holds
for every genus $g\geq 6$.
\end{proposition}


\textbf{Proof. }For the first proof we make an induction on $g$. 
Consider the boundary maps: 
$$
\xi _{a,A}: \Mgn{a,A\cup \left\{ s\right\} }\times
\Mgn{g-a,A^C\cup \left\{ t\right\} }
\rightarrow  \Mgn{g,P}, 
$$
on varying $(a,A)$ in such a way that $a\geq 3,g-a\geq 3$. 
Consider the composition of the induced
pull-back map with the projection on $H^2\otimes H^2$:
$$
g_{a,A}: \hmgn{4}{g,P} \rightarrow 
\hmgn{2}{a,A\cup \left\{ s\right\} }
\otimes 
\hmgn{2}{g-a,A^C\cup \left\{ t\right\}}. 
$$

We need a few remarks:

\begin{itemize}
\item  Under the above hypotheses on genera, there are no relation among
tautological classes in 
$\hmgn{2}{a,A\cup \left\{s\right\} } \otimes 
 \hmgn{2}{g-a,A^C\cup\left\{ t\right\} } $. 

\item  Every class of the standard basis in 
$\hmgn{2}{a,A\cup \left\{ s\right\} }
 \otimes 
 \hmgn{2}{g-a,A^C\cup \left\{ t\right\} }$
(by the standard basis we mean the one described in  \cite{AC}), 
appears, with the suitable sign, as a summand in the pull-back of at most 
one tautological class of 
$ \hmgn{4}{g,P}$, with the exception of 
$-\psi _s \otimes \psi _t$, which is a summand both of 
$\xi_{a,A}^{*}\left( \psi |\delta _{a,A}\right) $ and 
$\xi_{a,A}^{*}\left(\delta _{a,A}|\psi \right) $.
This is a combinatorial remark which follows
from the description of pull-backs of section \ref{pullbacks}. In particular,
one should look at the description of the operations on graphs denoted by 
$f_{s,t}$ and $j_{s,t}$.

\item  Almost every essential tautological class $\alpha $ in 
$\bgp{4}{g,P}$  
satisfies $g_{a,A}\left( \alpha \right) \neq 0$
for at least one $(a,A)$ satisfying the hypotheses. This is also
a combinatorial remark, and it is based on the relative position
of boundary cycles in $\Mgn{g,P}$.    
The exceptions are: 
\begin{eqnarray*}
\kappa _2 , \psi _x \text{ for every }x\in P , \delta _{E\left(b,B\right) } , \\
\delta _{G\left( c,C,d,D\right) }\text{, if }c+d\leq 2 .
\end{eqnarray*}
\end{itemize}

\noindent Suppose there is a relation among essential tautological classes in 
$\hmgn{4}{g,P}$. Applying all the maps $g_{a,A}$, 
one obtains that many coefficients have to vanish. The relation
should then be: 
$$
c\kappa _2+\sum_{x\in P}c_x \psi^2 _x
+\sum c_{b,B}\delta _{E\left( b,B\right) }
+\sum_{c+d\leq 2}c_{c,C,d,D}\delta _{G\left(c,C,d,D\right) }=0 
$$
We pull it back with the map 
$$
\xi ^{*}:\hmgn{4}{g,P}\rightarrow \hmgn{4}{g-1,P\cup \left\{ q,r\right\} } , 
$$
and get 
\begin{eqnarray*}
c\kappa _2  +\sum_{x\in P}c_x\psi^2 _x  
+\sum c_{b,B}\left( \delta _{E\left(b,B\right) }
+\delta _{E\left( b-1,B\cup \left\{ q,r\right\} \right)}+...\right) && \\
+\sum_{c+d\leq 2}c_{c,C,d,D}
\left( \delta _{G\left( c-1,C\cup \left\{q,r\right\} ,d,D\right) }
+\delta _{G\left( c,C,d-1,D\cup \left\{ q,r\right\}\right) }
+\delta _{G\left( c,C,d,D\right) }\right) &=&0
\end{eqnarray*}

\noindent By induction hypothesis, the coefficients $c,c_x,c_{b,B}$ all have to
vanish. Every type $G$ class appears at most once as a summand in the image
of at type $G$ class. If we call ``critical'' the classes corresponding to
graphs $G\left( 0,A,0,B\right) $, i.e. the possibly unessential ones, we
observe that every non-critical class has at least one non-critical summand
in its pull-back. On the other hand, if  we extend the ordering of $P$ to an 
ordering for $P\cup \left\{ q,r\right\} $ imposing $\left\{ q,r\right\} $ to be
 the last two elements, then a basis of critical classes maps to a set of linearly
independent critical classes. Thus, the coefficients $c_{c,C,d,D}$ vanish.

\begin{flushright}
$\square $
\end{flushright}

The main tool used in the second proof is the map: 
$$
\xi ^{*}:\hmgn{4}{g,P}\rightarrow \hmgn{4}{g-1,P\cup \left\{ q,r\right\} }. 
$$

\noindent The combinatorics of tautological classes and pull-back formulas becomes
rather intricate, but nevertheless it suggests a partition of $\bgp{4}{g,P}$, 
corresponding to any given partition of $P$, which, inductively, turns out
to give a direct sum decomposition of the tautological group.

\begin{definition}
\begin{enumerate}
\item {\bf Pure boundary classes of type $E$  and $F$} \\ 
are essential pure boundary classes corresponding to graphs 
$F$ and $E\left(a,A\right) $. \\
They generate the subspace $\mathbf{W}_{EF}$ of $\taut{4}{g,P}$.

\item {\bf Pure boundary classes of type $H$ and $G$} \\
are essential pure boundary classes corresponding to graphs 
$H\left( a,A\right) $ and $G\left( a,A,b,B\right)$. \\
They generate the subspace $\mathbf{W}_{GH}$ of $\taut{4}{g,P}$.

\item {\bf $\Psi$- mixed classes} \\
are essential mixed boundary classes 
$\psi |\delta _{irr}$ and $\psi|\delta _{a,A}$, 
generating $\mathbf{W}_\Psi .$

\item {\bf $\Psi_I$-mixed classes} \\
are essential mixed boundary classes $\psi _i\delta _{irr}$ 
and $\psi_i\delta _{a,A}$, with $i\in I \cap A$, 
generating $\mathbf{W}_{\Psi I}.$

\item  {\bf $K$- mixed classes} \\
are essential mixed boundary classes 
$\kappa _1\delta _{irr}$ and $\kappa |\delta _{a,A}$,
generating $\mathbf{W}_K.$

\item  { \bf Mumford $K$ classes} \\
are essential classes 
$\left\{ 
\begin{array}{c}
\kappa _1^2,\kappa _2\text{, for }g\geq 6              \\ 
\kappa _1^2\text{, for }g=5\text{ and }g=4,P=\emptyset  \\ 
\emptyset \text{, for }g=4,P\neq \emptyset \text{, and }g\leq 3
\end{array}
\right. $, and generate $\mathbf{K}$.

\item  {\bf Mumford $\Psi_I$ classes } \\
are essential  classes $\left\{ 
\begin{array}{c}
\kappa _1\psi _i,\psi _i^2,\psi _i\psi _j\text{, for }g\geq 4 \\ 
\psi _i^2,\psi _i\psi _j\text{, for }g=3 \\ 
\emptyset \text{, for }g\leq 2
\end{array}
\right. $, with $i,j\in I$, and generate $\mathbf{\Psi }_I$.

\item  {\bf Mumford $\Psi_{IJ}$ classes} \\
are essential classes $\left\{ 
\begin{array}{c}
\psi _i\psi _j\text{, for }g\geq 3 \\ 
\emptyset \text{, for }g\leq 2
\end{array}
\right. $, with $i\in I$, $j\in J$, and generate $\mathbf{\Psi }_{IJ}$.
\end{enumerate}
\end{definition}

\begin{proposition}
\label{piudisette}Suppose that Theorem \ref{free} holds for $g=6$. Then it
holds for every genus $g \geq 6 $.
\end{proposition}

\textbf{Proof. }
Let $O=\left\{ q,r\right\} $, so that $P\cup \left\{q,r\right\} =P\cup O$.
Following formulas of section \ref{pullbacks}, we describe how the above
subspaces of $\taut{4}{g,P}$ behave with respect to the map 
$$
\xi ^{*}:\hmgn{4}{g,P} \rightarrow \hmgn{4}{g-1,P\cup \left\{ q,r\right\}} . 
$$
We write down the behavior for genus $g\geq 4$. 
When no confusion will arise, we will denote by the same letter the subspaces 
of the same type in $\hmgn{4}{g,P}$ 
and $\hmgn{4}{g-1,P\cup \left\{ q,r\right\} }$.

\begin{eqnarray*}
K          &\rightarrow & K \text{, for }g\geq 7         \\
\Psi _P    &\rightarrow & \Psi_P \text{, for }g\geq 5    \\
W_K        &\rightarrow & W_K+W_\Psi +W_{\Psi P} +W_{EF} 
                          +W_{GH} +W_{\Psi O} +\Psi_O     \\
W_\Psi     &\rightarrow & W_\Psi +W_{GH}+\Psi _O          \\
W_{\Psi P} &\rightarrow & W_{\Psi P}+W_{GH}+\Psi _{OP}    \\
W_{EF}     &\rightarrow & W_{EF}+W_{GH}+W_{\Psi O}        \\
W_{GH}     &\rightarrow & W_{GH}+W_{\Psi O} .
\end{eqnarray*}
We prove the Proposition by induction on $g$. Suppose that 
$$
\taut{4}{g-1,P\cup \left\{ q,r\right\} }
=W_{EF} \oplus W_{GH}\oplus W_\Psi \oplus W_{\Psi P}\oplus 
 W_{\Psi O} \oplus W_K\oplus \Psi_P\oplus \Psi _O\oplus \Psi _{OP} 
$$
and that every summand is freely generated by essential tautological
classes. We write down in block form the matrix of the map 
$$
\xi ^{*}:\hmgn{4}{g,P} \rightarrow \hmgn{4}{g-1,P\cup \left\{ q,r\right\}} . 
$$ 
$$
\begin{tabular}{|l|l|l|l|l|l|l|l|l|l|l|}
\hline
& $K$        & $\Psi _P$ & $W_K$ & $W_\Psi $ & $W_{\Psi P}$ & $W_{EF}$ &
 $W_{GH}$ &$W_{\Psi O}$ & $\Psi _O$ & $\Psi _{OP}$ \\ \hline
$K$          & $A$   & $0$ & $0$   & $0$   & $0$   & $0$   & $0$   & $0$   & $0$   
& $0$  \\ \hline
$\Psi _P$    & $0$   & $B$ & $0$   & $0$   & $0$   & $0$   & $0$   & $0$   & $0$   
& $0$  \\ \hline
$W_K$        & $0$   & $0$ & $C$   & $...$ & $...$ & $0$   & $...$ & $...$ & $...$
& $0$  \\ \hline
$W_\Psi $    & $0$   & $0$ & $0$   & $D$   & $0$   & $0$   & $...$ & $0$   & $...$ 
& $0$  \\ \hline
$W_{\Psi P}$ & $0$   & $0$ & $0$   & $0$   & $E$   & $0$   & $...$ & $0$   & $0$  
& $...$\\ \hline
$W_{EF}$     & $0$   & $0$ & $0$   & $0$   & $0$   & $F$   & $...$ & $...$ & $0$  
& $0$  \\ \hline
$W_{GH}$     & $0$   & $0$ & $0$   & $0$   & $0$   & $0$   & $G$   & $...$ & $0$   
& $0$  \\ \hline
\end{tabular}
. 
$$
We claim that the elements of $\bgp{4}{g,P}$ form a basis for 
$\taut{4}{g,P}$. Because of the form of the above matrix, it is sufficient to
check that every subset generating each subspace consists of independent
classes. For this, we look at blocks $A,...,G$, and check that each of them
has maximal rank, equal to the number of rows.
It is easy to see that $A$ and $B$ are both the identity matrix, whereas from

\begin{eqnarray*}
\kappa _1\delta _{irr}\rightarrow \kappa _1\delta _{irr}+... \text{, if} g \geq 5 \\
\kappa |\delta _{a,A} \rightarrow \left\{ 
             \begin{array}{c}
             \kappa |\delta _{a,A}+
             \kappa |\delta _{a-1,A\cup \left\{ q,r\right\} }
                    \text{, if }g-a\geq 1\text{, }a\geq 4   \\ 
             \kappa |\delta _{a-1,A\cup \left\{ q,r\right\} }
                    \text{, if }g=a\geq 4                   \\ 
             \kappa |\delta _{a,A}
             \text{, if }g-a\geq 1, a\leq 3                 \\ 
             0\text{, if }g=a\leq 3
             \end{array}
             \right. 
\end{eqnarray*}
we observe that $C$ has maximal rank for $g \geq 5$.

\noindent Similarly, $D$ and $E$ have maximal rank for $g \geq 3$, wheres $F$
has maximal rank for $g \geq 2$. 

\noindent As for the block $G$, from   
\begin{eqnarray*}
\delta _{H\left( a,A\right) }\rightarrow \left\{ 
            \begin{array}{c}
            \delta _{H\left( a,A\right) }
            +\delta _{H\left( a-1,A\cup \left\{ q,r\right\}\right) }+...
            \text{, if }g-1-a\geq 1\text{, }a\geq 2        \\ 
            \delta _{H\left( a,A\right) }
            +\frac 56\delta_{H\left( a-1,A\cup \left\{q,r\right\} \right) }+...
            \text{, if }g-1-a\geq 1\text{, }a=1            \\ 
            \delta _{H\left( a-1,A\cup \left\{ q,r\right\} \right) }+...
            \text{, if } g=a+1\geq 3                       \\ 
            \frac 56\delta _{H\left( a-1,A\cup \left\{ q,r\right\} \right) }+...
            \text{, if }g=a+1=2                            \\ 
            \delta _{H\left( a,A\right) }+...
            \text{, if }g-1-a\geq 1\text{, }a=0            \\ 
            0\text{, if }g=a+1=1
\end{array}
\right. \\
\delta _{G\left( a,A,b,B\right) }\rightarrow 
\delta _{G\left( a,A,b,B\right)}
+\delta _{G\left( a-1,A\cup \left\{ q,r\right\} ,b,B\right) }
+\delta_{G\left( a,A,b-1,B\cup \left\{ q,r\right\} \right) };
\end{eqnarray*}
we observe that type $H$ classes are
independent, and independent from type $G$ ones. 
For the type $G$ class, the argument used in the proof
of Proposition \ref{dim1} works in this case as well. One can
write the block $G$ in a triangular form, and see that it has maximal rank
for $g \geq 3$.

\begin{flushright}
$\square $
\end{flushright}

\begin{lemma}
\label{sei}Suppose that Theorem \ref{free} holds for $g=5$. Then it
holds for genus $g=6.$
\end{lemma}

\textbf{Proof. }The same proof of Proposition \ref{piudisette} can be repeated
to prove that $\bgp{4}{6,P} \backslash \left\{ \kappa _2\right\} $ is
a set of linearly independent classes. Thus, if a relation does exist, 
it should be of the form 
$$
\kappa _2+...=0; 
$$
since $\xi ^{*}\left( \kappa _2+...\right) =\kappa _2+...=0$, then the
relation should be a pull-back of the relation in 
$\hmgn{4}{5,0}$ (see section \ref{newrel}): 
$$
\kappa_2-\frac 1{180}\delta _F
+\frac{37}{180}\delta _{E\left( 1\right)}+...=0, 
$$
and hence it should be of the form 
$$
\kappa _2-\frac 1{180}\delta _F
+\frac{37}{180}\left( \delta _{E\left(1,q\right) }
+\delta _{E\left( 1,r\right) }\right) +...=0; 
$$
but one can easily observe that classes $\delta _F$ and 
$\delta _{E\left(1,q\right) }+\delta _{E\left( 1,r\right) }$ 
do only appear in the pull-back 
$\xi ^{*}\left( \delta _F\right) =\delta _F+\left( \delta _{E\left(1,q\right) }
+\delta _{E\left( 1,r\right) }\right) +...$, hence cannot have
different coefficients. This leads to a contradiction.

\begin{flushright}
$\square $
\end{flushright}

\begin{proposition}
\label{rel32} There is a unique new relation in $\Mgn{3,2}$, and it is the one described in section \ref{newrel}.
\end{proposition}


\textbf{Proof.} 

We know from \cite{Fa1} and \cite{Fa}  the relations arising in $\hmgn{4}{3,0}$ and $\hmgn{4}{3,1}$,
and further we know that a new relation does exist in $\hmgn{4}{3,2}$, involving
pure Mumford classes $\psi_a$, $\psi_b$, $\psi_a \psi_b$.
We need to prove that the relation has exactly the form described in section \ref{newrel},
and that no other relation appears.
We also recall that the group $\hmgn{4}{2,2} $ has been computed in \cite{G2}.

The relations in $\hmgn{4}{3,0}$ and $\hmgn{4}{3,1}$ can be all used to write classes
$\kappa_1^2$, $\kappa_2$, $\kappa_1 \psi_i$ in terms of other boundary classes,
when $|P| \geq 2$.

Therefore, a possible new relation in $\hmgn{4}{3, 2} $, can be written as follows: 
\begin{equation}
\sum_\Gamma c_\Gamma \delta _\Gamma 
+\sum_pc_{p\left( irr\right) }p|\delta_{irr}
+\sum_{p, (a,A)}c_{p\left( a,A\right) }p|\delta _{a,A}
+\sum c_i\psi_i^2+ c_{ab}\psi _a\psi _b=0.  \label{relazione}
\end{equation}

The first constraints on coefficients in (\ref{relazione})
are derived by writing  down explicitly the non-vanishing pull-backs of tautological
classes under the map 
$$
\Mgn{3,s}\rightarrow \Mgn{3,ab}, 
$$
which glues a fixed rational tail marked by $P \cup {t}$ 
by identifying $t$ and  $s$, and observing that the pull-back in $\hmgn{4}{3,s} $ 
of (\ref{relazione}) must be  a multiple of Faber's 
relation involving $\kappa _1\psi _s$ (see section \ref{newrel}).
They are:

\begin{tabular}{llll}
$c_F=-\frac 1{630}k$ & $c_{H\left( 2,\emptyset \right) }=\frac 17k$ 
& $c_{G\left( 1,P,1,\emptyset \right) }=\frac{16}{35}k$ 
& $c_{\psi \left(irr\right) }=-\frac 1{42}k$ \\ 
$c_{E\left( 1,P\right) }=-\frac 9{35}k$ 
& $c_{H\left( 1,\emptyset \right)}=\frac 4{105}k$ 
& $c_{G\left( 1,\emptyset ,1,P\right) }=-\frac 8{35}k$ 
& $c_{\psi \left( 3,\emptyset \right) }=5k$ \\ 
$c_{E\left( 0,P\right) }=\frac{13}{21}k$ 
& $c_{H\left( 1,P\right) }=-\frac2{105}k$ 
& $c_{G\left( 1,\emptyset ,2,\emptyset \right) }=\frac 57k$ 
& $c_{\psi \left( 2,\emptyset \right) }=-\frac{40}{21}k$ \\ 
& $c_{H\left( 0,\emptyset \right) }=\frac 4{63}k$ 
& $c_{G\left( 2,\emptyset,0,P\right) }=\frac{40}{21}k$
& $c_{\psi \left( 2,P\right) }=-\frac{16}{21}k$ \\ 
& $c_{H\left( 0,P\right) }=\frac{10}{63}k$ 
& $c_{G\left( 2,\emptyset ,1,\emptyset\right) } =-k$ 
& $c_{\kappa \left( 3,\emptyset \right) }=-k$.
\end{tabular}

To determine the coefficient of some classes of type $H$ and $G$ we also need to use
the map 
$$
\hmgn{4}{3,P} \rightarrow \hmgn{2}{2,s}\otimes \hmgn{2}{1,P\cup t}.
$$

We then know by \cite{Fa} and \cite{Fa5} that a new relation 
does actually exist, and
therefore we fix the value of the constant $k$ to be $1$.

We consider the following maps:
\begin{eqnarray*}
\hmgn{4}{3,ab} & \rightarrow & \hmgn{2}{2,s} \otimes \hmgn{2}{1,abt} \\
\hmgn{4}{3,ab} & \rightarrow & \hmgn{2}{2,as} \otimes \hmgn{2}{1,bt} \\
\hmgn{4}{3,ab} & \rightarrow & \hmgn{4}{2,as};
\end{eqnarray*}
the constraints on the coefficient derived by pulling back (\ref{relazione})
force all of them  to be the ones indicated in section \ref{newrel}.

\begin{flushright}
$\square $ 
\end{flushright}


\begin{lemma}
\label{due} Theorem \ref{free} holds for genus $g=2.$
\end{lemma}


\textbf{Proof. } The cases $n=0,1$ are well known (see \cite{Mu}); the cases 
$n=2,3$ are entirely described in \cite{G2} and \cite{BP}. Recall that a new
relation appears in $\hmgn{4}{2,3} $ (see section \ref{newrel}).

\noindent For every set 
$\left\{i,j,k\right\} \subset P$, 
only the relation pulled back from $\Mgn{2,\left\{ i,j,k\right\} }$ 
contains the summand: 
$$
\psi _i\delta _{2,P\backslash \left\{ j,k\right\} }
+\psi _j\delta_{2,P\backslash \left\{ i,k\right\} }
+\psi _k\delta_{2,P\backslash \left\{i,j\right\} } ;
$$
we fix an ordering on $P$, and use the relation in 
$\hmgn{4}{2,\left\{ i,j,k\right\} }$ to express 
$\psi _i\delta_{2,P\backslash \left\{ j,k\right\} }$, 
for $i<j,i<k,$ as linear combination  of other classes. 

Let $\mathcal{C}^{4}_{2,P}$ be the set obtained from 
the set of essential classes $\bgp{4}{2,P}$ after having
eliminated the relations arising in degree $4$, that is, after
having removed all pure Mumford classes, and the classes
$\psi _i\delta_{2,P\backslash \left\{ j,k\right\} }$, 
for $i<j,i<k.$
Observe that the definition of $\mathcal{C}^{4}_{2,P}$
depends on the choice of an ordering on $P$.

If  $n=4$, there is no new relation among essential tautological
classes; we postpone the proof of this fact.  If $n\geq 5$, let $F_{2,P}^4$
be the free vector space generated by 
classes in $\mathcal{C}^{4}_{2,P}$.
One can define every pull-back map on $F_{2,P}^4$, following formulas in
section \ref{pullbacks}. Our claim is that the map 
$$
f=\left\{ f_{ij}^{*}\right\}: F_{2,P}^4 \longrightarrow 
\oplus_{\left\{ i,j\right\} \subset P}
F_{2,P\backslash \left\{i,j\right\} \cup \left\{ s\right\} }^4
$$
is injective for $|P|\geq 5$. This implies, by induction, that no new relation
among tautological classes can appear for $n\geq 5$: any new one should map
to zero with $f$.

We use a decomposition of $F_{2,P}^4$ similar to the one described at the
beginning of this section.

\begin{itemize}
\item  $W_F$ is generated by $\delta _F,$

\item  $W_E$ is generated by classes $\delta _{E\left( 1,A\right) },$

\item  $W_{H\left( 0\right) }$ is generated by classes $\delta _{H\left(
0,A\right) },$

\item  $W_{H\left( 1\right) }$ is generated by classes $\delta _{H\left(
1,A\right) },$

\item  $W_{G\left( 2,0\right) }$ is generated by classes $\delta _{G\left(
2,A,0,B\right) },$

\item  $W_{G\left( 0,2\right) }$ is generated by classes $\delta _{G\left(
0,A,2,B\right) },$

\item  $W_{G\left( 1,1\right) }$ is generated by classes $\delta _{G\left(
1,A,1,B\right) },$

\item  $W_{G\left( 1,0\right) }$ is generated by classes $\delta _{G\left(
1,A,0,B\right) },$

\item  $W_\psi $ is generated by classes $\psi |\delta _{2,A},$

\item  $W_{\psi _I}$ is generated by classes $\psi _i\delta _{2,A},$ with $
i\in I.$
\end{itemize}

In the space 
$\oplus_{\left\{ i,j\right\} \subset P}F_{2,P\backslash\left\{ i,j\right\} 
\cup \left\{ s\right\} }^4$, we denote by $W_X=\oplus_{ij}$ $W_X^{ij}$ 
the direct sum of subspaces 
$W_X^{ij}\subset F_{2,P\backslash \left\{ i,j\right\} \cup \left\{ s\right\} }^4$. 
The matrix of the map $f$ can be written in triangular block form 
(we omit all zeroes):

\noindent \begin{tabular}{|l|l|l|l|l|l|l|l|l|l|l|}
\hline
& $W_F$ & $W_E$ & $W_{H\left( 0\right) }$ & $W_{G\left( 1,0\right) }$ & $
W_{G\left( 0,2\right) }\oplus W_{\psi _S}$ & $W_{G\left( 2,0\right) }$ & $
W_{H\left( 1\right) }$ & $W_{G\left( 1,1\right) }$ & $W_\psi $ & $W_{\psi _P}
$ \\ \hline
$W_F$ & $A$ &  &  &  &  &  &  &  &  &  \\ \hline
$W_E$ &  & $B$ &  &  &  &  &  &  &  &  \\ \hline
$W_{H\left( 0\right) }$ &  &  & $C$ &  &  &  &  &  &  &  \\ \hline
$W_{G\left( 1,0\right) }$ &  &  &  & $D$ &  &  &  &  &  &  \\ \hline
$W_{G\left( 0,2\right) }$ &  &  &  &  & $E$ &  &  &  &  &  \\ \hline
$W_{G\left( 2,0\right) }$ &  &  &  &  &  & $F$ &  &  &  &  \\ \hline
$W_{H\left( 1\right) }$ & $...$ & $...$ & $...$ &  &  &  & $G$ &  &  &  \\ 
\hline
$W_{G\left( 1,1\right) }$ &  &  & $...$ & $...$ &  &  &  & $H$ &  &  \\ 
\hline
$W_\psi $ & $...$ & $...$ & $...$ & $...$ &  & $...$ &  &  & $I$ &  \\ \hline
$W_{\psi _P}$ & $...$ & $...$ & $...$ & $...$ & $...$ & $...$ & $...$ & $...$
& $...$ & $L$ \\ \hline
\end{tabular}

\smallskip

We just need to check that the blocks on the diagonal have maximal rank.
This is completely trivial for the blocks $A,B,C,D,E$. We check block 
$G$, and observe that blocks $H$ and $I$ present a very similar
combinatorics. $G$ is of the form 
$$
\left( 
\begin{array}{ccccc}
G^{12} & G^{13} & ... & G^{ij} & ...
\end{array}
\right) ,
$$
where $G^{ij}$ is a block of the matrix of the map $f_{ij}^{*}.$ We can
write $G^{ij}$ as 
$$
\begin{tabular}{|l|l|l|}
\hline
& $\delta _{H\left( 1,B\cup \left\{ s\right\} \right) }$ & 
$\delta _{H\left(1,B\right) }$ \\ \hline
$\delta _{H\left( 1,A\right) },\left\{ i,j\right\} \subset A$ & $Id$ & $0$
\\ \hline
$\delta _{H\left( 1,A\right) },\left\{ i,j\right\} \subsetneq A^C$ & $0$ &
 $Id$ \\ \hline
$\delta _{H\left( 1,P\backslash \left\{ i,j\right\} \right) }$ & $0$ & $...$
\\ \hline
$\delta _{H\left( 1,A\right) },|\left\{ i,j\right\} \cap A|=1$ & $0$ & $0$
\\ \hline
\end{tabular}
$$
We consider the matrix $G^{\prime }$ obtained removing the second column of
blocks from each $G^{ij}$, except for the columns corresponding to 
$\delta _{H\left( 1,\emptyset \right) },\delta _{H\left( 1,x\right) }.$
Finally, we can extract such a triangular matrix 
$$
\begin{tabular}{|l|l|l|}
\hline
& $\delta _{H\left( 1,B\right) },|B|\leq 1$ & 
$\delta _{H\left( 1,B\cup\left\{ s\right\} \right) }$ \\ \hline
$\delta _{H\left( 1,A\right) },|A|\leq 1$ & $Id$ & $0$ \\ \hline
$\delta _{H\left( 1,A\right) },|A|\geq 2$ & $...$ & $Id$ \\ \hline
\end{tabular}
.
$$
Observe that we just need the weaker assumption $|P|\geq 4$.

As for the block $L$, observe that any essential class  maps to 
essential  classes, except for  
$$
\psi _i\delta _{2,P\backslash \left\{ j,k\right\} }
\stackrel{f_{jk}^{*}}{\rightarrow}
 -\psi _i\psi _s
=-\sum_{|C^C|\geq 3}\psi _i\delta_{2,C}
 -\sum_{x<i}\psi _i\delta _{2,P\backslash \left\{ x,s\right\}
}+\sum_{x>i}\psi _x\delta _{2,P\backslash \left\{ i,s\right\} };
$$
but this doesn't prevent us from extracting a non-degenerate matrix 
$$
\begin{tabular}{|l|l|l|}
\hline
& $\psi _i\delta _{2,B},|B|\leq 1$ & 
$\psi _i\delta _{2,B\cup \left\{s\right\} }$ \\ \hline
$\psi _i\delta _{2,A},|A|\leq 1$ & $Id$ & $0$ \\ \hline
$\psi _i\delta _{2,A},|A|\geq 2$ & $...$ & $Id$ \\ \hline
\end{tabular}.
$$

With the same argument, one can write a sub-block of $F$ of the form 
\smallskip
$$
\begin{tabular}{|l|l|l|}
\hline
& $\delta _{G\left( 2,A,0,D\cup \left\{ s\right\} \right) },
   \delta _{G\left( 2,A,0,D\right) }$ 
& $\delta _{G\left( 2,C\cup \left\{ s\right\} ,0,D\right) }
$ \\ \hline
$\delta _{G\left( 2,A,0,B\right) },|A|\leq 1$ & $K$ & $0$ \\ \hline
$\delta _{G\left( 2,A,0,B\right) },|A|\geq 2$ & $...$ & $Id$ \\ \hline
\end{tabular}
$$

\smallskip

The set $P\backslash \left\{ i,j\right\} \cup \left\{ s\right\} $ inherits
an ordering from $P$, assuming $s$ to be the last point; therefore the
second column of blocks gives no problem. As for the matrix $K$, write it in
sub-blocks $K_A$, where $K_A$ involves classes $\delta _{G\left( 2,A,0,B\right)}$. 
These classes are all obtained pushing forward from 
$\hmgn{2}{0,A^C\cup \left\{ z\right\} } $, 
and so are the relations among them in $F_{2,P}^4$. 
The combinatorics of the map
corresponding to the block $K_A$ is then exactly the same of the map 
$$
\hmgn{2}{0,A^C\cup \left\{ z\right\} }
\rightarrow \oplus_{\left\{ i,j\right\} \subset A^c}
\hmgn{2}{0,A^C\backslash \left\{ i,j\right\} \cup \left\{ z,s\right\}} 
$$
which will  be proved in lemma \ref{inj0h2}  to be injective for 
$|A^C|\geq 4$. Therefore each $K_A$, and consequently $K$, has maximal rank.

As for the case $n=4$, we first prove by using the
pull-back map
\begin{equation*}
\hmgn{4}{2,\left\{ i,j,k,l\right\} }
\rightarrow 
\hmgn{2}{2,\left\{ i,s\right\} }
\otimes \hmgn{2}{0,\left\{ j,k,l,t\right\}} ;
\end{equation*}
that in a possible new relation, the coefficients of $\psi$-mixed classes and of  classes of type 
$G$ vanish.

We now restrict the map $f$ to the free vector space generated by the classes
with non vanishing coefficient in a possible new relation in
$\taut{4}{2,4}$.
By the same arguments used for the general case, the new map $f$ is injective,
and the proof of our Lemma is complete.

\begin{flushright}
$\square $
\end{flushright}


\begin{lemma}
\label{inj0h2}For $|P|\geq 5$, the map 
$$
\hmgn{2}{0,P} \rightarrow 
\oplus_{\left\{x,y\right\} \subset P\backslash \left\{ h\right\} }
\hmgn{2}{0,P\backslash \left\{ x,y\right\} \cup \left\{ s\right\}} 
$$
is injective.
\end{lemma}


\textbf{Proof.} The case $|P|=5$ is trivial. 

We can consider 
$$
\phi _A^{*}:\hmgn{2}{0,P} \rightarrow
H^2 ( \Mgn{0,A\cup \left\{ s\right\} }\times \Mgn{0,A^C\cup \left\{ t\right\} } ) 
$$
as the sum of the two maps 
$$
f_A^{*}:\hmgn{2}{0,P} \rightarrow 
\hmgn{2}{0,A\cup \left\{ s\right\} },
$$
$$
f_{A^C}^{*}:\hmgn{2}{0,P} \rightarrow
\hmgn{2}{0,A^c\cup \left\{ t\right\} },
$$
where the two maps are the pull-back of the map that glues any fixed rational
tail to the extra marked point. For any such $A$, there exist 
$\left\{x,y\right\} \subset P$ such that 
$A\subset P\backslash \left\{ x,y\right\} $. 
For a suitable choice of the rational tail to glue, we can write a
commutative diagram 
$$\xy
\xymatrix{
\Mgn{0,P\backslash \left\{ x,y\right\} \cup \left\{h\right\} } 
\ar[r]^{f_{P\backslash \left\{ x,y\right\} }}
& \Mgn{0,P} \\ 
\Mgn{0,A\cup \left\{ s\right\} } \ar[u] \ar[ur]^{f_A} & }
\endxy
$$
so that from the induced diagram on $H^2$ we read: 
$\ker f_{P\backslash\left\{ x,y\right\} }^{*}\subset \ker f_A^{*}$. 
Therefore, by proposition 2.8 in \cite{AC},
$$
\left( \cap _{\left\{ x,y\right\} \subset P}
\ker f_{P\backslash \left\{x,y\right\} }^{*}\right) 
\subset \left( \cap _{A\subset P,}\ker f_A^{*}\right) =0 .
$$
The statement is proved by induction on $|P|$ . Suppose  that 
$x\in \left(\cap _{\left\{ x,y\right\} \subset 
P\backslash \left\{ h\right\} }
\ker f_{P\backslash \left\{ x,y\right\} }^{*}\right) $, but there exist 
$k\in P\backslash \left\{ h\right\} $, such that 
$y\in f_{P\backslash \left\{h,k\right\} }^{*}\left( x\right) \neq 0$. 
By the commutativity of 

$$\xy
\xymatrix{
\hmgn{2}{0,P} \ar[r]^{f_{P\backslash\left\{ x,y\right\} }^{*}} \ar[d]^{f_{P\backslash \left\{ h,k \right\} }^{*}} 
& \hmgn{2}{0,P\backslash \left\{ x,y\right\} \cup \left\{ t\right\} }  \ar[d]\\
\hmgn{2}{0,P\backslash \left\{ h,k\right\} \cup\left\{ u\right\} } \ar[r]
& \hmgn{2}{0,P\backslash \left\{ h,k,x,y\right\} \cup \left\{ u,v\right\}} }
\endxy
$$

we see that $y\in \left( \cap _{\left\{ x,y\right\} 
\subset P\backslash\left\{ h,k\right\} }
\ker f_{P\backslash \left\{ x,y\right\} }^{*}\right) $,
hence by induction hypothesis, $y=0$, and we are done.

\begin{flushright}
$\square $
\end{flushright}


\textbf{Proof }of Theorem \ref{free}. The induction on the genus starts with Lemma
\ref{due}; then one can perform the next few steps by arguments similar to the one
used  in Lemma \ref{sei}, and get the result for genus up to $6$.
The procedure is then completed  with Proposition \ref{piudisette}.

\begin{flushright}
$\square $
\end{flushright}

\section{A conjecture on higher degree tautological relations}

At this point it is natural to
formulate a conjecture which is suggested by the proof
of Proposition \ref{dim1}. 

This conjecture agrees with Harer's and Ivanov's stability theorems 
(see \cite{Ha} and \cite{Iv}), and with Faber's
results and conjectures concerning  the tautological ring of the open
part $\mgno$ (see \cite{Fa5}).

\begin{conjecture}\label{conj}
There are no relations between essential tautological classes in
$\hmgn{2k}{g,P} $ whenever $g \geq 3k$.
\end{conjecture}

We justify our conjecture. We first need to extend some definitions.
A tautological class of degree $2k$ is a push-forward of a degree $2l$ 
Mumford class from a codimension $k-l$ boundary component; a degree $2k$
 class is unessential 
if it can be eliminated by means of a relation among tautological
classes arising in degree  $<2k$. 

Then we need to build new pull-back formulas, but we can give conjectural ones 
starting from the ones we proved for degree $4$. In particular, 
we claim that they preserve the tautological group.

Under the above hypotheses, there are plenty of boundary components 
$ \prod \Mgn{\Gamma_i} $ in
$\Mgn{g,P} $ such that
$$ H^{2k}( \prod \Mgn{\Gamma_i} ) = 
\oplus_{\sum i = 2k}( \otimes \hmgn{i}{g_i, P_i} ) 
$$
contains at least one summand
$\otimes \hmgn{2j}{g_j, P_j}$
with $g_j \geq  3j $, and $g_j < g$.

Write a generic linear combination of tautological
classes in $\hmgn{2k}{g,P} $, and suppose it is equal to $0$; 
by pulling back these relation to the above components
we can prove that many coefficients do vanish: in fact,
inductively, there are no relations among essential classes in these
summands of the cohomology.
We also conjecture  that the pull-back maps
in higher degree still satisfy the property that each class is generically a 
summand in the pull-back of at most one class.

It is then hard to believe that a new relation holds among the few
classes whose coefficient has not yet been showed to be zero.

\section{Generators of the cohomology group\label{generat}}

\begin{theorem}
\label{main} $\hmgnq{4}{g,P}$ is
generated by tautological classes for all $g \geq 8$.
\end{theorem}

\textbf{Proof}. We are following Edidin's scheme of Proof (\cite{Ed}). 

In the proof of this Proposition we plan to give
an upper bound for the dimension of
the cohomology group, and then to use the knowledge of the tautological group
and of the homology of the mapping class group to prove that,
this bound is achieved.

Let $n=3g-3+|P|$ be the complex dimension
of $\Mgn{g,P}$. 
We write a part of the exact homology sequence of the
pair $ (\Mgn{g,P} ,\Mgn{g,P}\backslash \Mgno{g,P})$ : 
$$
...\rightarrow H_{2n-4}\left( \Mgn{g,P}\backslash \Mgno{g,P}\right) 
   \stackrel{j_*}{\rightarrow} H_{2n-4}\left( \Mgn{g,P} \right)
   \rightarrow H_{2n-4}\left( \Mgn{g,P},\Mgn{g,P}\backslash \Mgno{g,P}\right) 
   \rightarrow ... 
$$
hence, using Poincar\'e duality for smooth orbifolds:  
\begin{eqnarray*}
\dim \hmgn{4}{g,P}&=& \dim H_{2n-4}\left( \Mgn{g,P} \right) \\
&\leq&   \dim j_* H_{2n-4}\left( \Mgn{g,P}\backslash \Mgno{g,P}\right)
       + \dim H^4\left( \Mgno{g,P}\right)
\end{eqnarray*}

\noindent We refer to the description of the stratified structure 
of $\Mgn{g,P}$ which has been explained in section \ref{notdef}.
For any stable graph $\Gamma$, we further denote by $\Delta_{\Gamma}^0$
the open stratum $\xi_{\Gamma} ( \Mgno{\Gamma} )$.

Let 
$$
\partial \Mgno{g,P}=\Mgn{g,P} \backslash \Mgno{g,P}.
$$ 
We recall that $\partial \Mgno{g,P}= \cup _i\Delta _{\Gamma _i}$, 
where the $\Delta _{\Gamma _i}$'s are the
codimension $1$ boundary components.

We denote by 
$ \partial\partial \Mgno{g,P}$ 
the union of the codimension two boundary
components, and write the homology exact sequence for the pair
$\left( \partial \Mgno{g,P},\partial \partial \Mgno{g,P}
\right) : $
$$
...\rightarrow H_{2n-4}\left( \partial \partial \Mgno{g,P} \right) 
    \stackrel{i_{*}}{\rightarrow} H_{2n-4}\left( \partial \Mgno{g,P} \right)
   \rightarrow H_{2n-4}\left( \partial \Mgno{g,P},\partial \partial \Mgno{g,P} \right) 
   \rightarrow ... 
$$

\noindent Let us look at the relative term. By Lefschetz Theorem (\cite{Sp})
we have:
$$
H_{2n-4}\left( \partial \Mgno{g,P},\partial \partial \Mgno{g,P} \right) 
\simeq
H^2 \left(\partial \Mgno{g,P} \backslash \partial \partial \Mgno{g,P} \right).
$$

\noindent The space
$$
\partial \Mgno{g,P} \backslash \partial \partial \Mgno{g,P} 
$$
consists of the disjoint union of the interior parts of the codimension
$1$ boundary components, the $\Delta_{i}^0$'s.

\noindent We have a precise description of these $\Delta_{i}^0$'s 
as quotients of moduli spaces of smooth curves:
$$
\partial \Mgno{g,P} \backslash \partial \partial \Mgno{g,P} 
\simeq 
\sqcup_{a,A} 
(\Mgno{a,A \cup \{ s \}} \times \Mgno{g-a,A^c \cup \{ t \}}) / Aut \Gamma_{a,A}
\sqcup \Mgno{g-1,P \cup \{ qr \}}/ Aut \Gamma_{irr}
$$

\noindent The rational cohomology of such quotients satisfies:
$$
H^k\left( \Mgno{\Gamma_i}/ Aut \Gamma_i , \mathbb{Q}\right) \cong 
H^k\left( \Mgno{\Gamma_i}\mathbb{Q} \right) ^{Aut \Gamma_i}
$$
where we denote by $H^k\left( \Mgno{\Gamma_i}\right) ^{Aut \Gamma_i}$ the
invariants with respect to the induced $Aut \Gamma_i$ action on the cohomology.
In the case $k=2$, these invariants can be precisely described. 
The cohomology group
$ H^2\left( \Mgno{\Gamma_i}\right)$ is generated by Mumford classes of
degree $2$. The class $\kappa_1$ is fixed by the automorphism
group of any graph, whereas the $\psi_i$ classes, for $i$
a special point ,
are permuted by the group action in the obvious way. 

We then get
 
\begin{eqnarray*}
H^2 (\partial \Mgno{g,P} \backslash \partial \partial \Mgno{g,P} ) & \simeq &
\oplus_{a,A} H^2 (\Mgno{a,A \cup \{ s \}} \times \Mgno{g-a,A^c \cup \{ t \}})
                   ^{Aut \Gamma_{a,A}}
\oplus  H^2(\Mgno{g-1,P \cup \{ q,r \}})^{Aut \Gamma_{irr}}.
\end{eqnarray*}

\noindent At this point, the bound for the dimension of the cohomology group is:
\begin{eqnarray*}
\dim \hmgn{4}{g,P} &\leq & 
 + \sum_{a,A} \dim H^2 (\Mgno{a,A \cup \{ s \}} \times \Mgno{g-a,A^c \cup \{ t \}})
                   ^{Aut \Gamma_{a,A}} \\
&& + \dim  H^2(\Mgno{g-1,P \cup \{ q,r \}})^{Aut \Gamma_{irr}}  \nonumber 
+\dim H^4 (\Mgno{g,P}) 
+\dim i_*j_* H_{2n-4}\left( \partial \partial \Mgno{g,P} \right)
\end{eqnarray*}

The space $\partial \partial \Mgno{g,P}$ is the union of the codimension
two boundary components, which we will call $\Theta_i$'s. Their
complex dimension is $n-2$. An easy application of the
Maier-Vietoris exact sequence, shows that the obvious map
$$
k: \sqcup_i \Theta_i \longrightarrow \cup_i \Theta_i 
= \partial \partial \Mgno{g,P}
$$
from the disjoint union into the union of these components
induces the following isomorphism in homology:
$$
 \oplus_i H_{2n-4} (\Theta_i) \simeq
       H_{2n-4} (\partial \partial \Mgno{g,P}).
$$
Observe that for dimension reasons, $ \dim H_{2n-4} (\Theta_i)=1$.

\noindent We claim that 
$$
\dim i_*j_* H_{2n-4}\left( \partial \partial \Mgno{g,P} \right) \leq r
$$
where $r$ equals  the number of essential pure boundary classes.
This number differs from the number of codimension two boundary
components because of the presence of Keel's relations in genus $0$.
These relations live in the
second homology group of $\Mgn{0,n}$.

\noindent The push-forward induced 
by the map
$$
\Mgn{0,A \cup \{s \}} \rightarrow \Mgn{g,P}
$$
determines homological equivalences among codimension $2$ boundary components
of $\Mgn{g,P}$.

Let 
$$
\phi: \sqcup_i \Theta_i \longrightarrow \Mgn{g,P}
$$ 
be the collection of the inclusion maps of
the codimension $2$ 
boundary components.
By what we said above, the image of the map
$$
\phi_* : H_{2n-4} ( \sqcup_i \Theta_i ) \longrightarrow H_4 (\Mgn{g,P})
$$
has dimension less or equal than $r$.
Since $ \phi = k \circ i \circ j $, and $k_*$ is an isomorphism,
this implies that 
$$
\dim i_*j_* H_{2n-4}\left( \partial \partial \Mgno{g,P} \right) \leq r.
$$

Our final bound is:
\begin{eqnarray} \label{disdim}
\dim \hmgn{4}{g,P} &\leq & 
 \sum_{a,A} \dim H^2 (\Mgno{a,A \cup \{ s \}} \times \Mgno{g-a,A^c \cup \{ t \}})
                   ^{Aut \Gamma_{a,A}} \nonumber \\
&&+ \dim  H^2(\Mgno{g-1,P \cup \{ qr \}})^{Aut \Gamma_{irr}}  \\
&&+\dim H^4 (\Mgno{g,P}) \nonumber
+r 
\end{eqnarray}

By Ivanov (\cite{Iv}), Harer (\cite{Ha}), and Loojenga's (\cite{Lo1})
stability theorems for the homology of the mapping class group, 
$H^4\left( \Mgno{g,P} \right)$ is freely generated by Mumford classes, for $g\geq
8$.
  
Instead of computing the dimension of all the
cohomology groups involved in (\ref{disdim}), 
we proceed more indirectly. We show that there is a bijection
between the following two sets. On one hand, the set 
 $\mathcal{B}^{4}_{g,P}$, on the other, the set whose elements are
the $r$ pure boundary classes in $\mathcal{B}^{4}_{g,P}$ and the vectors
belonging to the natural bases of the cohomology vector spaces appearing
on the right hand side of the above inequality (\ref{disdim}).
The upper bound for the dimension of the cohomology group is therefore
achieved, and consequentely the tautological classes generate the cohomology group.

The bijection directly follows from the definition of
essential tautological classes: 
\begin{itemize}
\item pure Mumford classes in 
$ \mathcal{B}^{4}_{g,P} $  
correspond  to a basis for 
$$
H^4( \Mgno{g,P}),
$$ 
\item mixed boundary classes in $\mathcal{B}^{4}_{g,P}$ correspond to a basis for 
$$
\oplus_{a,A} H^2 (\Mgno{a,A \cup \{ s \}} \times \Mgno{g-a,A^c \cup \{ t \}})
^{Aut \Gamma_{a,A}}
\oplus  H^2(\Mgno{g-1,P \cup \{ qr \}})^{Aut \Gamma_{irr}},
$$
\item pure boundary classes  in $\mathcal{B}^{4}_{g,P}$
 are exactly $r$. 
\end{itemize}

This completes the proof.

\begin{flushright}
$\square $
\end{flushright}

\bigskip

\noindent 
Department of Mathematics \\
California Institute of Technology \\ 
91125 Pasadena, CA  \\
polito@caltech.edu 

\end{document}